\documentclass[12pt]{amsart}
\usepackage{amsmath, amssymb,graphics,multirow}
\usepackage{hyperref}

\newtheorem{theorem}{Theorem}[section]
\newtheorem{lemma}[theorem]{Lemma}

\newtheorem{conjecture}[theorem]{Conjecture}

\theoremstyle{definition}
\newtheorem*{remark}{Remark}

\newcommand{\bbK}{{\mathbb K}}
\newcommand{\bbQ}{{\mathbb Q}}
\newcommand{\bbR}{{\mathbb R}}
\newcommand{\bbZ}{{\mathbb Z}}

\newcommand{\cD}{{\mathcal D}}
\newcommand{\cN}{{\mathcal N}}

\def\ord{\operatorname{ord}}
\newcommand{\sgn}{\operatorname{sgn}}

\numberwithin{equation}{section}

\begin{document}

\title[Lang's conjecture and sharp height estimates for $y^{2}=x^{3}+b$]{Lang's conjecture and sharp height estimates for the elliptic curves $y^{2}=x^{3}+b$}

\author{Paul Voutier}
\address{}
\curraddr{London, UK}
\email{paul.voutier@gmail.com}

\author{Minoru Yabuta}
\address{}
\curraddr{Senri High School, 17-1, 2 chome, Takanodai, Suita, Osaka, 565-0861, Japan}
\email{rinri216@msf.biglobe.ne.jp}


\dedicatory{}


\begin{abstract}
For $E_{b}: y^{2}=x^{3}+b$, we establish Lang's conjecture on a lower bound for
the canonical height of non-torsion points along with upper and lower bounds for
the difference between the canonical and logarithmic height. These results are
either best possible or within a small constant of the best possible lower bounds.
\end{abstract}

\maketitle

\section{Introduction}

The canonical height, $\widehat{h}$ (defined in Section~2), on an
elliptic curve $E$ defined over a number field $\bbK$ is a measure of the
arithmetic complexity of points on the curve. It has many desirable properties.
For example, it is a positive definite quadratic form on the lattice
$E(\bbK)/(\text{torsion})$, behaving well under the group law on $E(\bbK)$.
See \cite[Chapter~VIII]{Silv2} and
\cite[Chapter~9]{BG} for more information on this height.

There is another important, and closely related, height function defined for
points on elliptic curves, the absolute logarithmic height (also defined in
Section~2). It has a very simple definition which makes it very easy to compute.

In this paper, we provide sharp lower bounds for the canonical height as well
as bounding the difference between the heights for a well-known and important
family of elliptic curves, the Mordell curves defined by $E_{b}: y^{2}=x^{3}+b$
where $b$ is a sixth-power-free integer (i.e., quasi-minimal Weierstrass equations
for all $E_{b}/\bbQ$).

\subsection{Lower bounds}

Lang's Conjecture proposes a lower bound for the heights of non-torsion points
on a curve which varies with the curve.

\begin{conjecture}[Lang's Conjecture]
\label{conj:lang}
Let $E/\bbK$ be an elliptic curve with minimal discriminant $\cD_{E/\bbK}$. There
exist constants $C_{1}>0$ and $C_{2}$, depending only on $[\bbK:\bbQ]$, such that
for all nontorsion points $P \in E(\bbK)$ we have
\[
\widehat{h}(P) > C_{1} \log \left( \cN_{\bbK/\bbQ} \left( \cD_{E/\bbK} \right) \right) + C_{2}.
\]
\end{conjecture}

See \cite[p.~92]{Lang} along with the strengthened version in \cite[Conjecture VIII.9.9]{Silv2}.

Such lower bounds have applications to counting the number of integral points
on elliptic curves \cite{HS}, questions involving elliptic divisibility
sequences \cite{EIS, Everest1, VY1} and several other problems.

Silverman \cite[Section~4, Theorem]{Silv1} showed that Lang's conjecture holds
for any elliptic curve with $j$-invariant non-integral for at most $R$ places
of $\bbK$ (note that this includes our curves, $E_{b}$, since their $j$-invariant
is $0$), but with $C_{1}$ dependent on $\bbK$ and $R$. Gross and Silverman
\cite[Proposition~3(3)]{GS} proved an explicit version of this result from
which it follows that for non-torsion points, $P$, on $E_{b}$, we have
\[
\widehat{h}(P) > 3 \cdot 10^{-14} \log \left| \Delta \left( E_{b} \right) \right|.
\]

Also Hindry and Silverman \cite{HS} proved an explicit version of Lang's
conjecture whenever Szpiro's ratio, $\sigma_{E/\bbK}$, of $E/\bbK$ is known.
Hence Lang's conjecture follows from Szpiro's conjecture (or the $ABC$ conjecture).
Subsequently, David \cite{David} and Petsche \cite{Petsche} improved Hindry and
Silverman's result. It can be shown that $\sigma_{E_{b}/\bbQ}<5$, hence from
Petsche's Theorem~2, for example, a weaker result than the above follows with
$3 \cdot 10^{-14}$ replaced by $2 \cdot 10^{-22}$.

However, these results for $E_{b}/\bbQ$ all follow from more general results.
By focusing specifically on $E_{b}/\bbQ$, much better results can be obtained.

When $b$ is a nonzero integer that is sixth-power-free, Krir \cite[Proposition~3.1]{K}
showed that for any non-torsion point, $P$,
\[
\widehat{h}(P) > 10^{-3}\log |b| + 10^{-3}.
\]

In the special case of $b=-432m^{2}$ for a cube-free integer $m$, Jedrzejak
\cite{Jed} proved a sharper result, which was improved by Everest, Ingram and
Stevens \cite[Lemma~4.3]{EIS} and further improved very recently by Fujita and
Nara \cite[Proposition~2.5]{FN2}:
\[
\widehat{h}(P) > \frac{1}{18}\log |b|-1.1009.
\]
The coefficient of $\log|b|$ is correct in their result, but as we show below
in Theorem~\ref{thm:lang}(c), the constant should be $-(2/9)\log(2)-(1/4)\log(3)=-0.4286\ldots$.

Also if $b$ is a positive square-free integer, Fujita and Nara \cite[Proposition~4.3]{FN1}
showed that
\[
\widehat{h}(P) > \frac{1}{24}\log |b| - 0.073576,
\]
upon noting that their canonical height is twice ours (compare with our results
for this case in Theorem~\ref{thm:lang}(a) or \eqref{eq:lang1-pos} below).

We express the hypotheses of our theorem in terms of the {\it Tamagawa index
at $p$} for $p$, a prime. Letting $E_{0} \left( \bbQ_{p} \right)$ be the
connected component of the identity in $E \left( \bbQ_{p} \right)$,
the Tamagawa index, $c_{p}$, at $p$, is the order of the component group,
$E \left( \bbQ_{p} \right)/E_{0} \left( \bbQ_{p} \right)$, of $E$ at $p$.
See \cite{CPS} and \cite[Section~IV.9]{Silv6} for more details.

\begin{theorem}
\label{thm:lang}
Let $b$ be an integer which is sixth-power-free and let $P \in E_{b}(\bbQ)$
be a non-torsion point.

\noindent
{\rm (a)} If $c_{p}=1$ for all primes, $p>3$, then
\[
\widehat{h}(P) >
\left\{
	\begin{array}{ll}
		 \displaystyle\frac{1}{6}\log|b| - \log(2) - \frac{1}{2}\log(3) & \text{if $b<0$}\rule[-4.0mm]{0cm}{10.0mm} \\
		 \displaystyle\frac{1}{6}\log |b| - \frac{2}{3}\log(2) - \frac{3}{4}\log(3) - 0.006 & \text{if $b>0$.}\rule[-4.0mm]{0cm}{10.0mm}
	\end{array}
\right.
\]

\noindent
{\rm (b)} If $c_{p}|4$ for all primes, $p>3$ and $2|c_{p}$ for at
least one such prime, then
\[
\widehat{h}(P) >
\left\{
	\begin{array}{ll}
		 \displaystyle\frac{1}{24}\log |b| - \frac{1}{4}\log(2) - \frac{5}{48}\log(3) & \text{if $b<0$}\rule[-4.0mm]{0cm}{10.0mm} \\
		 \displaystyle\frac{1}{24}\log |b| - \frac{1}{6}\log(2) - \frac{1}{6}\log(3) - 0.002 & \text{if $b>0$.}\rule[-4.0mm]{0cm}{10.0mm}
	\end{array}
\right.
\]

\noindent
{\rm (c)} If $c_{p}|3$ for all primes, $p>3$ and $c_{p}=3$ for at
least one such prime, then
\[
\widehat{h}(P) >
\left\{
	\begin{array}{ll}
		 \displaystyle\frac{1}{18}\log |b| - \frac{2}{9}\log(2) - \frac{1}{4}\log(3) - 0.004 & \text{if $b<0$}\rule[-4.0mm]{0cm}{10.0mm} \\
		 \displaystyle\frac{1}{18}\log |b| - \frac{1}{3}\log(2) - \frac{1}{6}\log(3) - 0.004 & \text{if $b>0$.}\rule[-4.0mm]{0cm}{10.0mm}
	\end{array}
\right.
\]

\noindent
{\rm (d)} If $c_{p}|12$ for all primes, $p>3$, $2|c_{p}$ for at least one such
prime, $p$, and $3|c_{q}$ for at least one other such prime, $q$, then
\[
\widehat{h}(P) >
\left\{
	\begin{array}{ll}
		 \displaystyle\frac{1}{36} \log |b| - 0.2247 & \text{if $b<0$\rule[-4.0mm]{0cm}{10.0mm}} \\
		 \displaystyle\frac{1}{36} \log |b| - 0.2262 & \text{if $b>0$.\rule[-4.0mm]{0cm}{10.0mm}}
	\end{array}
\right.
\]
\end{theorem}

\begin{remark}
In the course of the proof of Theorem~\ref{thm:lang}, we establish the minimum
value of $\widehat{h}(P)$ for all possibilities of $b$ modulo powers of $2$
and $3$. As such bounds can be important for obtaining sharp results for other
problems (e.g., primitive divisor problems for elliptic divisibility sequences),
we refer the reader to these bounds in \eqref{eq:lang1-neg} and \eqref{eq:lang1-pos}
for part~(a), \eqref{eq:lang2-neg} and \eqref{eq:lang2-pos} for part~(b),
\eqref{eq:lang3-neg} and \eqref{eq:lang3-pos} (as well as \eqref{eq:lang3-bNeg-a}
and \eqref{eq:lang3-bNeg-b}) for part~(c) and
\eqref{eq:lang4-neg} and \eqref{eq:lang4-pos} for part~(d).

All that is required to apply these bounds is knowing the congruence classes of
$b$ modulo powers of $2$ and $3$, the reduction of $P$ (or $[2]P$ for part~(b))
at $2$ and $3$ and then referring to Tables~\ref{table:quant-2} and \ref{table:quant-3}.
\end{remark}

In \cite{VY2}, we were able to show that our results are best possible. See
Section~\ref{sect:sharp} for examples showing that Theorem~\ref{thm:lang}(a)
and (b) for $b<0$ are the best possible results and the other lower bounds are
within 0.006 of the best possible result. By ``best possible'', we mean that
the value for $C_{1}$ in Conjecture~\ref{conj:lang} is best possible and then,
fixing $C_{1}$, the value for $C_{2}$ is best possible.

The constants in part~(d) are not as ``nice'' as the ones in parts~(a)--(c) of
the theorem, but they do arise in a natural way in this setting. For example,
the best-possible constant for $b<0$ is $0.19155\ldots-(1/3)\log(2)-(1/6)\log(3)$
and $-0.19155\ldots$ is the minimal value of $(1/2)\log(c)-(1/12)\log \left( c^{5}-c^{2} \right)$
plus the sum in \eqref{eq:arch-hgt}, where $c$ is defined by $x(P)=c|b|^{1/3}$.
The log terms here again arise naturally in the proof of the theorem.

As in \cite{VY2}, our proof is based on the decomposition of the canonical height
as the sum of local height functions. However, there are differences in the
behaviour of the local height functions for the curves in each family. One of particular interest
to us is that the local archimedean height function here has an error term near
its critical point that is $O \left(\epsilon^{2} \right)$, whereas in \cite{VY2},
the analogous error term is $O(\epsilon)$. In general, it appears that the
archimedean height function for all elliptic curves behaves in one of these two
ways near critical points. Our work to understand this function better is
ongoing.

Lastly, note that while the formulation of our results is not in terms of
$\Delta \left( E_{b} \right)$, it is equivalent to such a formulation since
$\Delta \left( E_{b} \right)=-432b^{2}$.

\subsection{Difference of heights}

Our proof of our lower bound for the canonical height also allows us to prove
sharp bounds on the difference between the canonical height and
the logarithmic height of points on $E_{b}(\bbQ)$.

In Example~2.1 of \cite{Silv5}, Silverman showed that
\[
-\frac{1}{6}\log |b| - 1.576
\leq \frac{1}{2}h(P)-\widehat{h}(P)
\leq \frac{1}{6}\log |b| + 1.48
\]
and that the coefficients on $\log|b|$ are best possible.

Using a combination of Proposition~5.18(a) and Theorem~5.35(c) of \cite{S-Z},
one can obtain
\[
-\frac{1}{6}\log |b|-0.578 \leq \frac{1}{2}h(P)-\widehat{h}(P)
\leq \frac{1}{6}\log |b|+1.156.
\]

\begin{theorem}
\label{thm:hgt-diff}
Let $b$ be a nonzero integer and let $P \in E_{b}(\bbQ)$.

For $b<0$,
\[
-\frac{\log(3)}{4} - 0.005
< \frac{1}{2}h(P)-\widehat{h}(P)
< \frac{1}{6}\log |b| + \frac{\log(2)}{3} + \frac{\log(3)}{4},
\]
and for $b>0$,
\begin{eqnarray*}
&   & -\frac{1}{6}\log |b| - \frac{\log(2)}{3} - 0.007-0.076b^{-1/3} \\
& < & \frac{1}{2}h(P)-\widehat{h}(P) \\
& < & \frac{1}{6}\log |b| + \frac{\log(2)}{3} + \frac{\log(3)}{4} + 0.004.
\end{eqnarray*}

In all cases, we have
\[
-\frac{1}{6}\log |b| - 0.299 < \frac{1}{2}h(P)-\widehat{h}(P) < \frac{1}{6}\log |b| + 0.51.
\]
\end{theorem}

\begin{remark}
Only the upper bound when $b<0$ is best possible here. It appears that the terms
$-0.005$ in the lower bound for $b<0$, $-0.007$ in the lower bound for $b>0$ and
$0.004$ in the upper bound for $b>0$ are not required. Examples
demonstrating these claims are provided in Section~\ref{sect:sharp}.
\end{remark}

\begin{remark}
As with Theorem~\ref{thm:lang}, improved results can often be obtained for
specific congruence classes of $b$ modulo powers of $2$ and $3$, here by using
\eqref{eq:diff-neg} and \eqref{eq:diff-pos}.
\end{remark}

\section{Notation}

For what follows in the remainder of this paper, we will require some standard
notation (see \cite[Chapter 3]{Silv2}, for example).

Let $\bbK$ be a number field and let $E/\bbK$ be an elliptic curve given by the
Weierstrass equation
\[
E: y^{2} + a_{1}xy + a_{3}y = x^{3} + a_{2}x^{2} + a_{4}x + a_{6},
\]
with $a_{1},\ldots,a_{6} \in \bbK$.

Put
\begin{align*}
b_{2} & = a_{1}^{2}+4a_{2},  & b_{6} & = a_{3}^{2}+4a_{6}, \\
b_{4} & = 2a_{4}+a_{1}a_{3}, & b_{8} & = a_{1}^{2}a_{6}+4a_{2}a_{6}-a_{1}a_{3}a_{4}+a_{2}a_{3}^{2}-a_{4}^{2}.
\end{align*}

Then $E/\bbK$ is also given by $y^{2}=4x^{3}+b_{2}x^{2}+2b_{4}x+b_{6}$.

Furthermore, for any $m \in \bbZ$, $[m]:E/\bbK \rightarrow E/\bbK$ is the
multiplication-by-$m$ isogeny.

For a point $P \in E(\bbK)$, we define the canonical height of $P$ by 
\[
\widehat{h}(P)=\frac{1}{2}\lim_{n \to \infty}\frac{h \left( \left[ 2^{n} \right] (P) \right)}{4^{n}},
\]
with $h(P)=h(x(P))$, where $h(P)$ and $h(x(P))$ are the absolute logarithmic
heights of $P$ and $x(P)$, respectively (see Sections~VIII.6,7 and 9 of
\cite{Silv2}). Also recall that for $\bbQ$, $h(s/t)=\log\max \{|s|, |t|\}$
with $s/t$ in lowest terms is the absolute logarithmic height of $s/t$.

Let $M_{\bbK}$ be the set of valuations of $\bbK$ and for each $v \in M_{\bbK}$,
let $n_{v}$ be the local degree and let $\widehat{\lambda}_{v}(P):
E \left( \bbK_{v} \right) \backslash \{ O \} \rightarrow \bbR$ be the local
height function, where $\bbK_{v}$ is the completion of $\bbK$ at $v$. From
\cite[Theorem~VI.2.1]{Silv6}, we have the following decomposition of the
canonical height into local height functions
\[
\widehat{h}(P) = \sum_{v \in M_{\bbK}} n_{v} \widehat{\lambda}_{v}(P).
\]

For $\bbK=\bbQ$, the nonarchimedean valuations on $\bbK$ can be identified with
the set of rational primes. For a nonarchimedean valuation, $v$, we let
$q_{v}$ be the associated prime,
\[
v(x) = -\log |x|_{v} = \ord_{q_{v}}(x)\log \left( q_{v} \right)
\]
for $x \neq 0$ and $v(0)=+\infty$.

\begin{remark}
We refer the reader to \cite[Section~4]{CPS} and \cite[Remark~VIII.9.2]{Silv2}
for notes about the various normalisations of both the canonical and local height
functions. In what follows, our local height functions, $\widehat{\lambda}_{v}(P)$,
are those that \cite{CPS} denotes as $\lambda_{v}^{{\rm SilB}}(P)$, that is
as defined in Silverman's book \cite[Chapter~VI]{Silv6}. So, as stated in (11) of
\cite{CPS}, their $\lambda_{v} \left( P \right)$ equals
$2\widehat{\lambda}_{v}\left( P \right)+(1/6) \log \left| \Delta \left( E \right) \right|_{v}$
here.

Our canonical height also follows Silverman and is half that found in \cite{CPS}
as well as half that returned from the height function, ellheight, in PARI.
\end{remark}

\section{Archimedean Estimates}

\subsection{$b<0$}

\begin{lemma}
\label{lem:arch-b-neg}
Suppose $b \in \bbR$ is negative and let $P=(x(P),y(P)) \in E_{b}(\bbR)$ be a point
of infinite order.

\noindent
{\rm (a)} We have
\begin{equation}
\label{eq:archim-hgt-b-neg1}
\widehat{\lambda}_{\infty}(P)
> -\frac{1}{3}\log(2)=\frac{1}{6} \log |b| + \frac{1}{4} \log(3)
- \frac{1}{12} \log \left| \Delta \left( E_{b} \right) \right|.
\end{equation}

\noindent
{\rm (b)} Suppose that $x(P)=c|b|^{1/3}$ where $c>1$. Then
\begin{eqnarray}
\label{eq:archim-hgt-b-neg3a}
\hspace*{3.0mm} \widehat{\lambda}_{\infty}(P)
& > & \frac{1}{3} \log(c) - \frac{1}{36}\log(6912) -0.004 \\
& = & \frac{1}{6} \log |b| + \frac{1}{3} \log(c) + \frac{\log(108)}{18}
-0.004
- \frac{1}{12} \log \left| \Delta \left( E_{b} \right) \right| \nonumber
\end{eqnarray}
and
\begin{eqnarray}
\label{eq:archim-hgt-b-neg3b}
\widehat{\lambda}_{\infty}(P)
& > & \frac{1}{12} \log \left( \frac{c^{5}-c^{2}}{432} \right) + 0.1895 \\
& = & \frac{1}{6} \log |b| + \frac{1}{12} \log \left( c^{5}-c^{2} \right)
+ 0.1895 - \frac{1}{12} \log \left| \Delta \left( E_{b} \right) \right|. \nonumber
\end{eqnarray}

\noindent
{\rm (c)} We have
\[
-\frac{\log(3)}{4} - 0.005
< \left( \frac{1}{2} \log \max \left\{ 1, |x(P)| \right\}
-\frac{1}{12} \log \left| \Delta \left( E_{b} \right) \right| \right)
- \widehat{\lambda}_{\infty}(P)
< 0.
\]
\end{lemma}

\begin{remark}
We have expressed the bounds in parts~(a) and (b) both with and without the
$(1/12)\log \left| \Delta \left( E_{b} \right) \right|$ term. The former
expression will be used in the proof of our theorems, while the latter is of
interest as it demonstrates that these results are actually independent of $b$.

All of the bounds are either best possible or within at most 0.005 of the
best possible results.

The actual dependence of $\widehat{\lambda}_{\infty}(P)$ on $c$
as $c \rightarrow +\infty$ is $(1/2)\log(c)$, but the lower bounds in \eqref{eq:archim-hgt-b-neg3a}
and \eqref{eq:archim-hgt-b-neg3b} allow us to obtain near best possible lower
bounds for the canonical height.
\end{remark}

\begin{proof}
We will estimate the archimedean contribution to the canonical height by using
Tate's series (see \cite{Tate2} as well as the presentation in \cite{Silv3}). Let
\[
t(P)=1/x(P) \quad \text{and} \quad
z(P) = 1-b_{4}t(P)^{2}-2b_{6}t(P)^{3}-b_{8}t(P)^{4},
\]
for a point $P=(x(P), y(P)) \in E(\bbR)$. Then the archimedean local height of
$P \in E(\bbR)$ is given by the series
\begin{equation}
\label{eq:arch-hgt}
\widehat{\lambda}_{\infty}(P)
= \frac{1}{2} \log |x(P)|
+ \frac{1}{8} \sum_{k=0}^{\infty} 4^{-k} \log \left| z \left( \left[ 2^{k} \right](P) \right) \right|
- \frac{1}{12} \log \left| \Delta(E) \right|,
\end{equation}
provided $x \left( \left[ 2^{k} \right](P) \right) \neq 0$ for all $k \neq 0$.

Here we have $b_{2}=b_{4}=b_{8}=0$ and $b_{6}=4b$, so $t(P)=1/x(P)$ and
$z(P) = 1-8bt(P)^{3}$.

Since $b<0$, for any $Q \in E_{b}(\bbR)$, $x(Q) \geq |b|^{1/3}$.
Hence $1 \leq z(Q) \leq 9$. In particular, $1 \leq z \left( \left[ 2^{k} \right](P) \right) \leq 9$.

(a) Applying this inequality for $k \geq 1$ and the definition of $z(P)$ to
\eqref{eq:arch-hgt}, we obtain
\[
0 \leq \widehat{\lambda}_{\infty}(P)
- \left( \frac{1}{8} \log \left( x(P)^{4}-8bx(P) \right) - \frac{1}{12} \log \left| \Delta \left( E_{b} \right) \right| \right)
\leq \frac{1}{12} \log(3).
\]

Since $x(P) \geq |b|^{1/3}$, $x(P)^{4}-8bx(P) \geq 9|b|^{4/3}$. Hence
\begin{eqnarray*}
\widehat{\lambda}_{\infty}(P)
& \geq & \frac{1}{8} \log \left( x(P)^{4}-8bx(P) \right) - \frac{1}{12} \log \left| \Delta \left( E_{b} \right) \right| \\
& \geq & \frac{1}{6} \log |b| + \frac{1}{4}\log(3) - \frac{1}{12} \log \left| \Delta \left( E_{b} \right) \right|
\end{eqnarray*}

(b) We first consider \eqref{eq:archim-hgt-b-neg3a}.
From \eqref{eq:arch-hgt}, our expression for $x(P)$ and our lower bound
for $z \left( \left[ 2^{k} \right](P) \right)$, we have
\[
\widehat{\lambda}_{\infty}(P)
\geq (1/6)\log |b| + (1/2)\log(c)
+ \frac{1}{8} \sum_{k=0}^{2} 4^{-k} \log \left| z \left( \left[ 2^{k} \right](P) \right) \right|
- \frac{1}{12} \log \left| \Delta(E) \right|,
\]
so we now proceed to bound from below
\begin{equation}
\label{eq:lb-3.3}
(1/6)\log(c)
+ \frac{1}{128} \log \left( z^{16} \left( P \right) z^{4} \left( [2](P) \right) z \left( [4](P) \right) \right).
\end{equation}

The derivative of this quantity is a rational function of $c$ whose numerator
is of degree $63$ with leading coefficient $1$ and whose denominator is of
degree $64$ with leading coefficient $6$. The numerator has only one root with
$c \geq 1$, which is $1.71216\ldots$, while the denominator has no roots with
$c \geq 1$. Therefore, the minimum value of \eqref{eq:lb-3.3} is
$(1/18)\log(108)-0.00372\ldots$, which occurs at $c=1.71216\ldots$.

We now consider \eqref{eq:archim-hgt-b-neg3b}. Proceeding as in the proof for
\eqref{eq:archim-hgt-b-neg3a}, we bound from below
\[
(1/2)\log(c) - (1/12) \log \left( c^{5}-c^{2} \right)
+ \frac{1}{128} \log \left( z^{16} \left( P \right) z^{4} \left( [2](P) \right) z \left( [4](P) \right) \right).
\]

The derivative of this quantity is a rational function of $c$ whose numerator
is of degree $66$ with leading coefficient $1$ and whose denominator is of
degree $67$ with leading coefficient $12$. The numerator has only one root with
$c \geq 1$, which is $4.21378\ldots$, while the denominator has no roots with
$c>1$. Therefore, the minimum value occurs at $4.21378\ldots$ and is
$0.1895\ldots$.

(c) We estimate
\[
\sum_{k=0}^{\infty} 4^{-k} \log \left| z \left( \left[ 2^{k} \right](P) \right) \right|.
\]

We will proceed in a similar way to the proof of part~(b). We group adjacent
triples of summands together and consider
$z^{16} \left( P \right) z^{4} \left( [2](P) \right) z \left( [4](P) \right)$,
which is a rational function of $c$ whose numerator and denominator
are both of degree $63$.

Neither the numerator nor the denominator of the derivative of this rational
function has a root with $c \geq 1$. So this function is
decreasing for $c \geq 1$. Hence it takes its maximum value,
$3^{32}$, at $c=1$ and its minimum value is $1$, which is approached from
above as $c \rightarrow +\infty$.

Therefore
\begin{eqnarray*}
0 &   <  & \frac{1}{8}\sum_{k=0}^{\infty} 4^{-k}\log \left( z \left( \left[ 2^{k} \right](P) \right) \right) \\
  & \leq & \frac{1}{8}\sum_{k=0}^{\infty}\frac{\log \left( 3^{32} \right)}{16 \cdot 64^{k}}
=\frac{16}{63}\log(3) = \frac{\log(3)}{4}+0.0043\ldots,
\end{eqnarray*}
so part~(c) follows from \eqref{eq:arch-hgt} and since
$\max \left\{ 1, |x(P)| \right\}=|x(P)|$ here.
\end{proof}

\subsection{$b>0$}

We use Tate's series here too. However for $b>0$, $E_{b}(\bbR)$ includes the point
$\left( 0,b^{1/2} \right)$, which causes a problem since we require $x \left( \left[ 2^{k} \right](P) \right)$
to be bounded away from $0$ to ensure that Tate's series converges. To get around
this, we use an idea of Silverman's (see \cite[p.~340]{Silv3}) and translate
the curve to the right using $x'=x+2b^{1/3}$, noting that $\widehat{\lambda}_{\infty}$
is fixed under such translations. In this way, we obtain the elliptic curve
\[
E_{b}': y^{2}=x^{3}-6 b^{1/3}x^{2}+12b^{2/3}x-7b
\]
and every point, $P'(x,y)$, in $E_{b}'(\bbR)$ satisfies $x\left( P' \right) \geq b^{1/3}$.
Here we have $b_{2}=-24b^{1/3}$, $b_{4}=24b^{2/3}$,
$b_{6}=-28b$ and $b_{8}=24b^{4/3}$. Hence
\[
t\left( P' \right)=1/x\left( P' \right) \quad \text{and} \quad
z\left( P' \right)=1-24b^{2/3}t\left( P' \right)^{2}+56bt\left( P' \right)^{3}-24b^{4/3}t\left( P' \right)^{4}.
\]

We could take the same approach here as in the proof of Lemma~3.4 of \cite{VY2}.
However there is a significant complication. Whereas in \cite{VY2},
for $x\left( P' \right)=(1+\epsilon)\sqrt{a}$, we have
\[
\widehat{\lambda}_{\infty}\left( P' \right)- \left\{ (1/4)\log(a)-(1/12)\log \left| \Delta \left( E_{a} \right) \right| \right\}
=O(\epsilon),
\]
here we find that for $x\left( P' \right)=2(1+\epsilon)b^{1/3}$,
\[
\widehat{\lambda}_{\infty}\left( P' \right)- \left \{ (1/6)\log(b)+(1/3)\log(2)-(1/12)\log \left| \Delta \left( E_{b} \right) \right| \right\}
=O \left( \epsilon^{2} \right),
\]
so we would need to proceed much more carefully. We have done so in an earlier
version of this paper. Here we take a more direct, although more computational,
approach and instead work with the actual expressions in the first terms of
Tate's series. The cost of this is a small additional constant term.

Despite this change in approach from earlier versions, the following result
may still be of interest to readers.

\begin{lemma}
\label{lem:bnds-b-pos-sharp}
Suppose $b \in \bbR$ is positive and $P' \in E_{b}'(\bbR)$ where $x\left( P' \right)=2(1+\epsilon)b^{1/3}$.

\noindent
{\rm (a)} For $-0.1745 \leq \epsilon \leq 0.6$,
\[
\log \left| x\left( P' \right)^{4}z\left( P' \right) \right|
\geq \log \left( 8b^{4/3} \right) - 2\epsilon - 2\epsilon^{2} + \frac{16}{3}\epsilon^{3}
 + 9.7\epsilon^{4},
\]

\noindent
{\rm (b)} Suppose $k$ is a positive integer and that $-0.379 \leq (-2)^{k-1}\epsilon \leq 1.044$.
Then
\[
\left| \frac{x \left( \left[ 2^{k} \right]\left( P' \right) \right)}{2b^{1/3}}
- \left\{ 1 + (-2)^{k}\epsilon
+ \left( (-2)^{4k} - (-2)^{k} \right) \epsilon^{4} \right\} \right|
\leq 2 \cdot 2^{7k}|\epsilon|^{7}.
\]

\noindent
{\rm (c)} If $k$ is a positive integer and $-1.0 \leq (-2)^{k} \epsilon \leq 0.36$, then
\begin{align*}
& \log \left( z \left( \left[ 2^{k} \right]\left( P' \right) \right) \right) \geq \\
& -\log(2) -6(-2)^{k} \epsilon + 4(-2)^{3k} \epsilon^{3}
+ \left( 9(-2)^{4k} + 6(-2)^{k} \right) \epsilon^{4}
+ (144/5) (-2)^{5k}\epsilon^{5}.
\end{align*}
\end{lemma}

\begin{proof}
These are parts~(a), (b) and (d) of Lemma~3.2 in\\
\verb+http://arxiv.org/abs/1305.6560v1+.
\end{proof}

\begin{lemma}
\label{lem:bnds-b-pos-1}
Suppose $b \in \bbR$ is positive and $P' \in E_{b}'(\bbR)$ be a point of
infinite order. If $K$ is a non-negative integer, then
\begin{equation}
\label{eq:rem-pos}
-1.8 \cdot 4^{-K}
< \sum_{k=K}^{\infty} 4^{-k} \log \left( z \left( \left[ 2^{k} \right]\left( P' \right) \right) \right)
< 2.24 \cdot 4^{-K}.
\end{equation}
\end{lemma}

\begin{remark}
These bounds are close to best possible. The best possible constants appear to
be $-1.7835\ldots$ when $x\left( P' \right)$ is near $2.9399\ldots b^{1/3}$ and
$\log(9)=2.197\ldots$ as $x\left( P' \right)$ approaches $b^{1/3}$.
\end{remark}

\begin{proof}
We write $x\left( P' \right)=2(1+\epsilon)b^{1/3}$, noting that $\epsilon \geq -0.5$.
We will group adjacent triples of summands together and consider
$$
f_{16}(\epsilon)=z^{16} \left( P' \right) z^{4} \left( [2]\left( P' \right) \right) z \left( [4]\left( P' \right) \right),
$$
which is a rational function of $\epsilon$ whose numerator and denominator
are both of degree $64$.

For $\epsilon \geq -0.5$, the numerator of the derivative of $f_{16}(\epsilon)$
has just one root, $\epsilon=0.41859\ldots$. The denominator of the derivative
is $262144(1+\epsilon)^{65}$, which is positive for $\epsilon \geq -0.5$. So
$f_{16}(\epsilon)$ is decreasing in the interval $-0.5 \leq \epsilon \leq 0.41859\ldots$
and increases towards $1$ for larger $\epsilon$. Hence $f_{16}(\epsilon)$ takes
its maximum value, $3^{32}$, at $\epsilon=-0.5$ and its minimum value,
$5.182\ldots 10^{-13}$, at $\epsilon=0.41859\ldots$.

Therefore
\begin{align*}
-1.8 \cdot 4^{-K}
& < 4^{-K}\sum_{k=0}^{\infty}\frac{\log \left( 5.182\ldots 10^{-13} \right)}
{16 \cdot 64^{k}} \\
& < \sum_{k=K}^{\infty} 4^{-k}\log \left( z \left( \left[ 2^{k} \right]\left( P' \right) \right) \right)
< 4^{-K}\sum_{k=0}^{\infty}\frac{\log \left( 3^{32} \right)}{16 \cdot 64^{k}}
<2.24 \cdot 4^{-K}.
\end{align*}
\end{proof}

\begin{lemma}
\label{lem:arch-b-pos}
Let $b \in \bbR$ be a positive real number and let $P \in E_{b}(\bbR)$ be a
point of infinite order.

\noindent
{\rm (a)} We have
\begin{eqnarray}
\label{eq:archim-hgt-b-pos1}
\widehat{\lambda}_{\infty}(P)
& > & -\frac{1}{4}\log(3)-0.006 \nonumber \\
& = & \frac{1}{6} \log |b| + \frac{1}{3} \log(2) - 0.006
      -\frac{1}{12} \log \left| \Delta \left( E_{b} \right) \right|.
\end{eqnarray}

\noindent
{\rm (b)} Suppose that $x(P)=c|b|^{1/3}$ with $c>-1$ and $c \neq 0$. Then
\begin{eqnarray}
\label{eq:archim-hgt-b-pos3a}
\widehat{\lambda}_{\infty}(P)
& > & \frac{1}{3}\log|c|-\frac{1}{3}\log(2)-0.004 \nonumber \\
& = & \frac{1}{6} \log |b| + \frac{1}{3} \log|c| + \frac{1}{4}\log(3) - 0.004
- \frac{1}{12} \log \left| \Delta \left( E_{b} \right) \right|,
\end{eqnarray}
and
\begin{eqnarray}
\label{eq:archim-hgt-b-pos3b}
\widehat{\lambda}_{\infty}(P)
& > & \frac{1}{12}\log \left( c^{5}+c^{2} \right)-\frac{1}{12}\log(432)+0.188 \nonumber \\
& = & \frac{1}{6} \log |b| + \frac{1}{12} \log \left( c^{5}+c^{2} \right)
+0.188 - \frac{1}{12} \log \left| \Delta \left( E_{b} \right) \right|.
\end{eqnarray}

\noindent
{\rm (c)} For $b \geq 2$,
\begin{align*}
&   -\frac{1}{6} \log(b) - \frac{\log(2)}{3}-0.007-0.076b^{-1/3} \\
& < \left( \frac{1}{2} \log \max \left\{ 1, \left| x(P) \right| \right\}
             -\frac{1}{12} \log \left| \Delta \left( E_{b} \right) \right|
        \right)
         - \widehat{\lambda}_{\infty}(P)
< 0.004.
\end{align*}
\end{lemma}

\begin{remark}
As in Lemma~\ref{lem:arch-b-neg}, we have expressed the bounds in parts~(a) and (b)
both with and without the $(1/12)\log \left| \Delta \left( E_{b} \right) \right|$
term.

All of the bounds are either best possible or within at most $0.007$
of the best possible results.

For \eqref{eq:archim-hgt-b-pos1}, our proof shows that the $-0.006$ term is not
required if $x\left( P' \right)=2(1+\epsilon)b^{1/3}$ with $\epsilon \leq -0.033$
or $\epsilon \geq 0.128$. It is only for $\epsilon$ approaching $0$ where a more
careful analysis is required to eliminate this term. Likewise, the small constant
terms in the other inequalities are only required for small intervals around the
location of the minimal value.

As in Lemma~\ref{lem:arch-b-neg}, the actual dependence of $\widehat{\lambda}_{\infty}(P)$ on $c$
as $c \rightarrow +\infty$ is $(1/2)\log(c)$, but the lower bounds in \eqref{eq:archim-hgt-b-pos3a}
and \eqref{eq:archim-hgt-b-pos3b} allow us to obtain near best possible lower
bounds for the canonical height.

Lastly, in part~(c), it appears that the correct lower order term is $O \left( b^{-2/3} \right)$,
not $O \left( b^{-1/3} \right)$.
\end{remark}

\begin{proof}
(a) Write $x\left( P' \right)=2(1+\epsilon)b^{1/3}$ where $\epsilon \geq -0.5$ (i.e., we use
the point where $\widehat{\lambda}_{\infty}\left( P' \right)$ takes its minimum value as the
centre).

We proceed in a similar way as in the proof of Lemma~\ref{lem:bnds-b-pos-1} and
consider $f_{64}(\epsilon)=x\left( P' \right)^{64}z\left( P' \right)^{16}z\left( [2]\left( P' \right) \right)^{4}z\left( [4]\left( P' \right) \right)/b^{64/3}$.
This is a polynomial in $\epsilon$ of degree $64$.

The derivative of $f_{64}(\epsilon)$ has two roots with $\epsilon \geq -0.5$,
one at $\epsilon=-0.48899\ldots$ and the second at $\epsilon=0.07196\ldots$.
The former is a local maximum, while the latter is a local minimum. We find
that
$$
\widehat{\lambda}_{\infty}\left( P' \right) - \frac{1}{6}\log |b| + \frac{1}{12}\log \left| \Delta(E) \right|
> \frac{\log \left( f_{64}(\epsilon) \right)}{8 \cdot 4^{2}}
- \frac{1.8}{8 \cdot 4^{3}} - \log(2)/3 =-0.0056\ldots
$$
(the second term coming from Lemma~\ref{lem:bnds-b-pos-1} with $K=3$), and the
desired inequality follows.

(b) We apply \eqref{eq:arch-hgt} with $P'$
rather than $P$, where $x\left( P' \right)=(c+2)|b|^{1/3}$ with $c \geq -1$ and
the lower bound in Lemma~\ref{lem:bnds-b-pos-1} with $K=3$, we have
\begin{eqnarray}
\label{eq:lem34b}
\widehat{\lambda}_{\infty}\left( P' \right)
& \geq & \frac{1}{6}\log |b| + \frac{1}{2}\log(c+2) - 1.8/512 - \frac{1}{12} \log \left| \Delta(E) \right| \nonumber \\
&      & + \frac{1}{128} \log \left( z^{16} \left( P' \right) z^{4} \left( [2]\left( P' \right) \right)
         z \left( [4]\left( P' \right) \right) \right).
\end{eqnarray}

For \eqref{eq:archim-hgt-b-pos3a}, we proceed similarly, using \eqref{eq:lem34b},
and bound from below
\[
(1/2)\log(c+2)-(1/3)\log(c)
+ \frac{1}{128} \log \left( z^{16} \left( P' \right) z^{4} \left( [2]\left( P' \right) \right)
z \left( [4]\left( P' \right) \right) \right).
\]

The derivative of this quantity is a rational function of $c$ whose numerator
is of degree $64$ with leading coefficient $1$ and whose denominator is of
degree $65$ with leading coefficient $6$. The numerator has only one root with
$c>0$, which is $1.71508\ldots$, while the denominator has no roots with
$c>0$. Therefore, the minimum value occurs at $c=1.71508\ldots$ and is
$(1/4)\log(3)-0.00029\ldots$, establishing \eqref{eq:archim-hgt-b-pos3a} for $c>0$.

For $-1 \leq c <0$, we proceed in the same way, but bound from below
\[
(1/2)\log(c+2)-(1/3)\log(-c)
+ \frac{1}{128} \log \left( z^{16} \left( P' \right) z^{4} \left( [2]\left( P' \right) \right)
z \left( [4]\left( P' \right) \right) \right).
\]

For \eqref{eq:archim-hgt-b-pos3b}, we again proceed similarly using \eqref{eq:lem34b},
and here bound from below
\[
(1/2)\log(c+2)-(1/12)\log \left( c^{5}+c^{2} \right)
+ \frac{1}{128} \log \left( z^{16} \left( P' \right)
z^{4} \left( [2]\left( P' \right) \right) z \left( [4]\left( P' \right) \right) \right).
\]

Using the derivative of this quantity, we find its minimum value occurs
at $c=3.6038\ldots$ and is $0.1880\ldots$.

(c) We handle separately the cases of $|x(P)|>1$ and $|x(P)| \leq 1$. In both
cases, we write $x(P)=cb^{1/3}$ and consider
\[
\log \left( \frac{x\left( P' \right)^{64}z\left( P' \right)^{16}z \left( [2]\left( P' \right) \right)^{4}z\left( [4]\left( P' \right) \right)}{\max \left\{ 1, \left| x(P) \right| \right\}^{64}} \right).
\]

So for $|x(P)| \leq 1$, we consider
$$
g_{64}(c)=x\left( P' \right)^{64}z\left( P' \right)^{16}z\left( [2]\left( P' \right) \right)^{4}z\left([4]\left( P' \right) \right).
$$

This is $b^{64/3}$ times a polynomial in $c$ of degree $64$ with
integer coefficients. For $c \geq -1$, $g_{64}'(c)$
has roots at $c=-0.9779\ldots$ and $c=0.1439\ldots$. The former is a local
maximum, while the latter is a local minimum. The value of $g_{64}(c)$ at this
local minimum is $5.28\ldots \cdot 10^{12}b^{64/3}$. Note that the local
maximum corresponds to $|x(P)| \leq 1$ only if $b<1.069\ldots$,
while the same holds for the local minimum for $b<335.59\ldots$. So for
$1.069\ldots < b < 335.59\ldots$, $g_{64}(c)$
decreases from $x(P)=-1$ to $x(P)=0.1439\ldots b^{1/3}$ (where its value is
greater than $1$) and then increases to $x(P)=1$. For $b>335.59\ldots$,
$g_{64}(c)$ is monotonically decreasing for $x(P) \in [-1,1]$.

Now we must examine the values of $g_{64}(c)$ corresponding to $x(P)=\pm 1$. At
$x(P)=-1$, we have $g_{64}(c)=2^{43}b^{64/3}+4\cdot 2^{43}b^{63/3}
-48 \cdot 2^{43}b^{61/3}+\cdots$ and we find that for $b \geq 18.5$,
$g_{64}(c) < 2^{45}b^{63/3}$, so
$$
\log \left( g_{64}(c) \right)
< \log \left( 2^{43}b^{64/3} \right) + 4b^{-1/3},
$$
since $\log(1+x)<x$ for $x>0$.

By considering the extrema found via
calculus, we find that
$$
\log \left( g_{64}(c) \right)
< \log \left( 2^{43}b^{64/3} \right) + 9.7b^{-1/3}
$$
for $2 \leq b \leq 18.5$ too.

At $x(P)=1$, we have $g_{64}(c)=2^{43}b^{64/3}-4\cdot 2^{43}b^{63/3}+ 48 \cdot 2^{43}b^{61/3}+\cdots$
and we proceed in the same way to show that
$\log \left( g_{64}(c) \right)
< \log \left( 2^{43}b^{64/3} \right) + 9.7b^{-1/3}$ for $b \geq 2$ here too. Furthermore,
$g_{64}(c)>1$ for such $b$ too. Therefore,
$$
0 < \log \left( g_{64}(c) \right) < (64/3)\log(b)+43\log(2)+9.7b^{-1/3}.
$$

Now we consider the case $|x(P)| \geq 1$ and
\[
\frac{x\left( P' \right)^{64}z\left( P' \right)^{16}z \left( [2]\left( P' \right) \right)^{4}z\left( [4]\left( P' \right) \right)}{\left| x(P) \right|^{64}},
\]
which is a rational function of $c$. The numerator of the derivative of this
rational function is of degree $63$, has $-128$ as its leading coefficient and
no roots with $c \geq -1$, while the denominator is of degree $65$, has $1$ as
its leading coefficient and only has a root at $c=0$. Therefore, the rational
function is increasing for $-1 \leq c <0$ and decreasing for $c>0$. Combining
this with the above results for $x(P)=\pm 1$ establishes
\begin{eqnarray*}
0 & < & \log \left( \frac{x\left( P' \right)^{64}z\left( P' \right)^{16}z \left( [2]\left( P' \right) \right)^{4}z\left( [4]\left( P' \right) \right)}{\max \left\{ 1, \left| x(P) \right| \right\}^{64}} \right) \\
  & < & (64/3)\log(b)+43\log(2)+9.7b^{-1/3},
\end{eqnarray*}
for $b \geq 2$.

Applying this to our expression for $\widehat{\lambda}_{\infty}\left( P' \right)$ in
\eqref{eq:arch-hgt}, we obtain
\begin{eqnarray*}
&   & (1/6)\log(b)-(43/128)\log(2)+0.076b^{-1/3} \\
& < & \left( \frac{1}{2} \log \max \left\{ 1, \left| x(P) \right| \right\}
      + \frac{1}{8} \sum_{k=3}^{\infty} 4^{-k} \log \left| z \left( \left[ 2^{k} \right]\left( P' \right) \right) \right|
      - \frac{1}{12} \log \left| \Delta \left( E_{b} \right) \right| \right) \nonumber \\
&   & - \widehat{\lambda}_{\infty}\left( P' \right) < 0.
\end{eqnarray*}

Using Lemma~\ref{lem:bnds-b-pos-1} with $K=3$, it follows that
\begin{eqnarray*}
&   & -\frac{1}{6} \log(b) - \frac{1}{3}\log(2)-0.007-0.076b^{-1/3} \\
& < & \left( \frac{1}{2} \log \max \left\{ 1, \left| x(P) \right| \right\}
             -\frac{1}{12} \log \left| \Delta \left( E_{b} \right) \right|
      \right)
      - \widehat{\lambda}_{\infty}\left( P' \right)
      < 0.004.
\end{eqnarray*}

Part~(c) now follows upon recalling that $\widehat{\lambda}_{\infty}$ is fixed
under translation, so the same inequalities also hold for
$\widehat{\lambda}_{\infty}(P)$.
\end{proof}

\section{Nonarchimedean estimates}

\subsection{Nonarchimedean estimates for $q_{v}>3$}

\begin{lemma}
\label{lem:curves-v-odd}
Let $v$ be a nonarchimedean valuation on $\bbQ$ associated with a prime
number, $q_{v}>3$, and let $b$ be an integer such that $q_{v}^{6} \nmid b$.
The Kodaira types and Tamagawa indices of $E_{b}$ at $v$ are as in
Table~$\ref{table:kodaira-odd-not-3}$.

\begin{table}[h]
\begin{tabular}{|c|c|c|}
\hline
$b$                                                        & {\rm Kodaira type}    & $c_{v}$ \\ \hline
$\ord_{q_{v}}(b)=0$                                        & $I_{0}$               & $1$     \\ \hline
$\ord_{q_{v}}(b)=1$                                        & $II$                  & $1$     \\ \hline
$\ord_{q_{v}}(b)=2$                                        &                       &         \\
$b/q_{v}^{2}$ {\rm a quadratic residue modulo} $q_{v}$     & $IV$                  & $3$     \\ \hline
$\ord_{q_{v}}(b)=2$                                        &                       &         \\
$b/q_{v}^{2}$ {\rm a quadratic non-residue modulo} $q_{v}$ & $IV$                  & $1$     \\ \hline
$\ord_{q_{v}}(b)=3$, $q_{v} \equiv 1 \bmod 6$              &                       &         \\
$b/q_{v}^{3}$ {\rm a cubic non-residue modulo} $q_{v}$     & $I_{0}^{*}$           & $1$     \\ \hline
$\ord_{q_{v}}(b)=3$, $q_{v} \equiv 5 \bmod 6$              & $I_{0}^{*}$           & $2$     \\ \hline
$\ord_{q_{v}}(b)=3$, $q_{v} \equiv 1 \bmod 6$              &                       &         \\
$b/q_{v}^{3}$ {\rm a cubic residue modulo} $q_{v}$         & $I_{0}^{*}$           & $4$     \\ \hline
$\ord_{q_{v}}(b)=4$                                        &                       &         \\
$b/q_{v}^{4}$ {\rm a quadratic residue modulo} $q_{v}$     & $IV^{*}$              & $3$     \\ \hline
$\ord_{q_{v}}(b)=4$                                        &                       &         \\
$b/q_{v}^{4}$ {\rm a quadratic non-residue modulo} $q_{v}$ & $IV^{*}$              & $1$     \\ \hline
$\ord_{q_{v}}(b)=5$                                        & $II^{*}$              & $1$     \\ \hline
\end{tabular}
\caption{$E_{b}$ reduction information for $q_{v}>3$}
\label{table:kodaira-odd-not-3}
\end{table}
\end{lemma}

\begin{proof}
We use Tate's algorithm with $K=\bbQ_{v}$ (using the steps and notation in
Silverman's presentation of Tate's algorithm in \cite[Section~IV.9]{Silv6}).

\vspace{1.0mm}

{\sc Step~1.} This step applies when $\ord_{q_{v}}\left( \Delta \left( E_{b} \right) \right) =0$.
Since $\Delta \left( E_{b} \right) = -432 b^{2}$ and $432=2^{4} \cdot 3^{3}$, the
reduction type is $I_{0}$ at $v$ when $\ord_{q_{v}}(b)=0$.

\vspace{1.0mm}

{\sc Step~2.} We have $\ord_{q_{v}}\left( \Delta \left( E_{b} \right) \right)>0$.
The singular point, $P=(x(P), y(P))$, is already at $(0,0)$ since
$\ord_{q_{v}}(2y(P)), \ord_{q_{v}} \left( 3x(P) \right)>0$ implies that $\ord_{q_{v}}(x(P))>0$
too, so no change of variables is needed. Therefore, $b_{2}=0$
and hence $\ord_{q_{v}} \left( b_{2} \right)>0$. Thus Step~2 does not apply.

\vspace{1.0mm}

{\sc Step~3.} Since $a_{6}=b$, if $\ord_{q_{v}} \left( b \right) = 1$,
then the reduction type is $II$.

\vspace{1.0mm}

{\sc Step~4.} We may now assume that $\ord_{q_{v}}(b) \geq 2$.
Note that $b_{6}=4b$ and $b_{8}=0$. Hence $\ord_{q_{v}} \left( b_{8} \right) \geq 3$
and so Step~4 cannot apply.

\vspace{1.0mm}

{\sc Step~5.} If $\ord_{q_{v}} \left( b \right) = 2$, then the reduction type
is $IV$. If $b/q_{v}^{2}$ is a quadratic residue modulo $q_{v}$, then $c_{v}=3$.
Otherwise, $c_{v}=1$.

\vspace{1.0mm}

{\sc Step~6.} We write $P(T)=T^{3}+b/q_{v}^{3}$, since $a_{2}=a_{4}=0$.
Its discriminant is $-27b^{2}/q_{v}^{6}$. If $\ord_{q_{v}} \left( b \right) = 3$,
then the discriminant is not zero modulo $q_{v}$ and the reduction type is $I_{0}^{*}$.

If $-b/q_{v}^{3}$ is a cubic residue modulo $q_{v}$, then $P(T)$ has at least one
root in $k$. Since $-1$ is always a cubic residue, this condition is equivalent
to $b/q_{v}^{3}$ being a cubic residue modulo $q_{v}$, so we will always consider
$b/q_{v}^{3}$ instead in what follows. Note that if $-3$ is a quadratic residue
modulo $q_{v}$ (that is, $q_{v} \equiv 1 \bmod 6$), then $P(T)$ has three roots
in $k$ and $c_{v}=4$, otherwise (that is, $q_{v} \equiv 5 \bmod 6$) it only has
one root in $k$ and $c_{v}=2$.

If $b/q_{v}^{3}$ is not a cubic residue modulo $q_{v}$, then $c_{v}=1$.
It is an easy consequence of Fermat's little theorem that this is only possible
for $q_{v} \equiv 1 \bmod 6$.

\vspace{1.0mm}

{\sc Step~7.} Here we assume that $P(T)$ has one simple root and one
double root. But the third roots of unity are distinct, since $q_{v}>3$, so this
is not possible.

\vspace{1.0mm}

{\sc Step~8.} Again, since the third roots of unity are distinct, this
can only occur if the triple root of $P(T)$ is zero. That is, $\ord_{q_{v}}(b)>3$.
So we consider the polynomial $Y^{2}-b/q_{v}^{4}$. It has distinct roots if and
only if $\ord_{q_{v}}(b)=4$. 

If $\ord_{q_{v}}(b)=4$ and $b/q_{v}^{4}$ is a quadratic residue modulo $q_{v}$,
then the reduction type is $IV^{*}$ and $c_{v}=3$.
If $\ord_{q_{v}}(b)=4$ and $b/q_{v}^{4}$ is a non-quadratic residue modulo $q_{v}$,
then the reduction type is $IV^{*}$ and $c_{v}=1$.

\vspace{1.0mm}

{\sc Step~9.} Since $a_{4}=0$, this step does not apply.

\vspace{1.0mm}

{\sc Step~10.} This is the last remaining case if $b$ is sixth-power-free.
Here the reduction type is $II^{*}$.

This completes the proof.
\end{proof}

\begin{lemma}
\label{lem:points-v-odd}
Let $v$ be a nonarchimedean valuation on $\bbQ$ associated with a prime
number, $q_{v}>3$, and let $b$ be an integer such that $q_{v}^{6} \nmid b$.

\noindent
{\rm (a)} $P \in E_{b} \left( \bbQ_{v} \right)$ has singular reduction if and
only if
$$
\ord_{q_{v}}(x(P)), \ord_{q_{v}}(y(P))>0.
$$

\noindent
{\rm (b)} For any $P \in E_{b}\left( \bbQ_{v} \right) \backslash \{ O \}$,
\begin{align}
\label{eq:vadic-not-3-hgt}
\widehat{\lambda}_{v}(P)
= & \frac{1}{2} \log \max\{ 1, \left| x(P) \right|_{v} \}
    - \frac{\log \left| \Delta \left( E_{b} \right) \right|_{v}}{12}
\nonumber \\
  & - \left\{
	\begin{array}{ll}
		(1/3) \log \left( q_{v} \right) & \text{if $\ord_{q_{v}} \left( x(P) \right)>0$, $\ord_{q_{v}}(b)=2$,} \\
		                                & \text{and $b/q_{v}^{2}$ is a quadratic residue modulo $q_{v}$} \\
		(1/2) \log \left( q_{v} \right) & \text{if $\ord_{q_{v}} \left( x(P) \right)>0$, $\ord_{q_{v}}(b)=3$,} \\
					                    & \text{and $b/q_{v}^{3}$ a cubic residue modulo $q_{v}$}\\
		(2/3) \log \left( q_{v} \right) & \text{if $\ord_{q_{v}} \left( x(P) \right)>0$, $\ord_{q_{v}}(b)=4$,} \\
					                    & \text{and $b/q_{v}^{4}$ a quadratic residue modulo $q_{v}$} \\
		0                               & \text{otherwise.}
	\end{array}
\right.
\end{align}
\end{lemma}

\begin{proof}
(a) We require $\ord_{q_{v}} \left( 3x(P)^2 \right)=2\ord_{q_{v}} \left( x(P) \right)>0$
and $\ord_{q_{v}} \left(2y(P) \right)=\ord_{q_{v}} \left(y(P) \right)>0$.

(b) This follows from our results in Lemma~\ref{lem:curves-v-odd} along with
Proposition~6 and the accompanying Table~2, as well as equation~(11), of
\cite{CPS}.
\end{proof}

\subsection{Nonarchimedean estimates for $q_{v}=3$}

\begin{lemma}
\label{lem:curves-v-3}
Let $b$ be an integer and suppose that $3^{6} \nmid b$. The Kodaira types and
Tamagawa indices of $E_{b}$ at $3$ are as in Table~$\ref{table:kodaira-3}$.
\end{lemma}

\begin{table}[h]
\begin{tabular}{|c|c|c|}
\hline
$b$                                    & {\rm Kodaira type} & $c_{3}$ \\ \hline
$b \equiv  2, 3, 4, 5, 6, 7 \bmod 9$   & $II$               & $1$     \\ \hline
$b \equiv  1, 8 \bmod 9$               & $III$              & $2$     \\ \hline
$b \equiv 9 \bmod 27$                  & $IV$               & $3$     \\ \hline
$b \equiv 18 \bmod 27$                 & $IV$               & $1$     \\ \hline
$b \equiv  54,  81, 108 \bmod 243$     & $IV^{*}$           & $3$     \\ \hline
$b \equiv 135, 162, 189 \bmod 243$     & $IV^{*}$           & $1$     \\ \hline
$b \equiv 27, 216 \bmod 243$           & $III^{*}$          & $2$     \\ \hline
$b \equiv 0 \bmod 243$                 & $II^{*}$           & $1$     \\ \hline
\end{tabular}
\caption{$E_{b}$ reduction information for $q_{v}=3$}
\label{table:kodaira-3}
\end{table}

\begin{proof}
As in the proof of the previous lemma, Tate's algorithm is used here. But we do
not provide all the details here. The conservative reader is referred to the
earlier version of this paper, \verb+http://arxiv.org/abs/1305.6560v2+,
which contains the full details. Also, since $\bbQ_{3}^{*}/\bbQ_{3}^{*6}$ is
finite group of small size, the reader can verify this lemma using an
implementation of Tate's algorithm like \verb+elllocalred+ in PARI.
\end{proof}

\begin{lemma}
\label{lem:points-v-3}
Let $b$ be an integer and suppose that $3^{6} \nmid b$.

\noindent
{\rm (a)} $P \in E_{b} \left( \bbQ_{3} \right)$ has singular reduction if and
only if $\ord_{3}(x(P)+b)>0$.

\noindent
{\rm (b)} For any $P \in E_{b}\left( \bbQ_{3} \right) \backslash \{ O \}$,
\begin{align}
\label{eq:3adic-hgt}
\widehat{\lambda}_{3}(P)
= & \frac{1}{2}\log \max \left\{ 1, \left| x(P) \right|_{3} \right\}
    - \frac{\log \left| \Delta \left( E_{b} \right) \right|_{3}}{12} \nonumber \\
  & - \left\{
	\begin{array}{ll}
		(1/4) \log(3) & \text{if $b \equiv 1,8         \bmod   9$ and $\ord_{3} \left( x(P)+b \right)>0$} \\
		(1/3) \log(3) & \text{if $b \equiv 9           \bmod  27$ and $\ord_{3} \left( x(P) \right)>0$} \\
		(2/3) \log(3) & \text{if $b \equiv 54, 81, 108 \bmod 243$ and $\ord_{3} \left( x(P) \right)>0$} \\
		(3/4) \log(3) & \text{if $b \equiv 27, 216     \bmod 243$ and $\ord_{3} \left( x(P) \right)>0$} \\
		0 & \text{otherwise.}
	\end{array}
\right.
\end{align}
\end{lemma}

\begin{proof}
(a) We require $\ord_{3} \left( 3x(P)^2 \right)>0$ and
$\ord_{3} \left( 2y(P) \right)=\ord_{3} \left( y(P) \right)>0$, so
$\ord_{3}(x(P)^{3}+b)>0$.

Writing $x(P)=x_{n}/x_{d}$, we have $x(P)^{3}+b= \left( x_{n}^{3}+bx_{d}^{3} \right)/x_{d}^{3}$.
Since $x^{3} \equiv x \bmod 3$ for all $x$, we see that
$x_{n}^{3}+bx_{d}^{3} \equiv x_{n}+bx_{d} \bmod 3$ and therefore
$\ord_{3} \left( x(P)^{3}+b \right)>0$ if and
only if $\ord_{3} \left( x(P)+b \right)>0$.

(b) This follows from Lemma~\ref{lem:curves-v-3}, along with Proposition~6
and equation~(11) of \cite{CPS}.
\end{proof}

\subsection{Nonarchimedean estimates for $q_{v}=2$}

\begin{lemma}
\label{lem:curves-v-2}
Let $b$ be an integer and suppose that $2^{6} \nmid b$. The Kodaira types and
Tamagawa indices of $E_{b}$ at $2$ are as in Table~$\ref{table:kodaira-2}$.
\end{lemma}

\begin{table}[h]
\begin{tabular}{|c|c|c|}
\hline
$b$                         & {\rm Kodaira type} & $c_{2}$ \\ \hline
$b \equiv  16    \bmod 64$  & --                 & $1$     \\ \hline
$b \equiv  2, 3  \bmod  4$  & $II$               & $1$     \\ \hline
$b \equiv  5     \bmod  8$  & $IV$               & $1$     \\ \hline
$b \equiv  1     \bmod  8$  & $IV$               & $3$     \\ \hline
$b \equiv  8, 12 \bmod 16$  & $I_{0}^{*}$        & $2$     \\ \hline
$b \equiv  4     \bmod 32$  & $IV^{*}$           & $3$     \\ \hline
$b \equiv 20     \bmod 32$  & $IV^{*}$           & $1$     \\ \hline
$b \equiv 32,48  \bmod 64$  & $II^{*}$           & $1$     \\ \hline
\end{tabular}
\caption{$E_{b}$ reduction information for $q_{v}=2$}
\label{table:kodaira-2}
\end{table}

\begin{remark}
For $b \equiv 16 \bmod 64$, the minimal model is given by $y^{2}+y=x^3+(b-16)/64$
and its Kodaira type is $I_{0}$.
\end{remark}

\begin{proof}
As above, we apply Tate's algorithm and refer the reader to either the earlier
version of this paper, \verb+http://arxiv.org/abs/1305.6560v2+ or the use of
an implementation of Tate's algorithm like \verb+elllocalred+ in PARI.

We only add that in Step~9 for $b \equiv 16 \bmod 64$, we need to apply the
translation $y=y'+4$, obtaining the curve $y^{2}+8y=x^{3}+b-16$. Neither Step~9
nor Step~10 apply, and in Step~11 we find that our Weierstrass equation is not
minimal and we obtain a new Weierstrass equation
\[
y'^{2}+y'=x'^{3}+(b-16)/64.
\]
The discriminant of this Weierstrass equation is $-27b^{2}/2^{8}$. Since this
is odd, the reduction type is $I_{0}$.
\end{proof}

\begin{lemma}
\label{lem:points-v-2}
Let $b$ be an integer and suppose that $2^{6} \nmid b$.

\noindent
{\rm (a)} $P \in E_{b} \left( \bbQ_{2} \right)$ has singular reduction if and
only if $\ord_{2}(x(P))>0$.

\noindent
{\rm (b)} For any $P \in E_{b}\left( \bbQ_{2} \right) \backslash \{ O \}$,
\begin{align}
\label{eq:2adic-hgt}
\widehat{\lambda}_{2}(P)
= & \frac{1}{2}\log \max \left\{ 1, \left| x(P) \right|_{2} \right\}
    - \frac{\log \left| \Delta \left( E_{b} \right) \right|_{2}}{12} \nonumber \\
  & - \left\{
	  \begin{array}{ll}
		(1/3) \log(2) & \text{if $b \equiv 1 \bmod 8$ and $\ord_{2}(x(P))>0$} \\
		(1/2) \log(2) & \text{if $b \equiv 8, 12 \bmod 16$ and $\ord_{2}(x(P))>0$} \\
		(2/3) \log(2) & \text{if $b \equiv 4 \bmod 32$ and $\ord_{2}(x(P))>0$} \\
		      \log(2) & \text{if $b \equiv 16 \bmod 64$ and $\ord_{2}(x(P))>0$} \\
		0 & \text{otherwise.}
	  \end{array}
\right.
\end{align}
\end{lemma}

\begin{remark}
Note that for $b \equiv 16 \bmod 64$, $E_{b}$ is not a minimal model. However,
since it arises in several cases, including the Mordell curve, $y^{2}=x^{3}-432m^{2}$
associated with cubic twists of the Fermat cubic, we include the result here.

Furthermore, this inclusion allows us to handle all $E_{b}$ by simply removing
any sixth-powers.
\end{remark}

\begin{proof}
(a) We require $\ord_{2} \left( 3x(P)^2 \right)
=2\ord_{2} \left( x(P) \right)>0$ and $\ord_{2} \left( 2y(P) \right)>0$. Since
$b \in \bbZ$ and $\ord_{2} \left( x(P) \right)>0$, $\ord_{2} \left( 2y(P) \right)>0$
always holds. Hence $\ord_{2} \left( x(P) \right)>0$ is a necessary and sufficient
condition.

(b) The proof is identical to that for part~(b) of Lemmas~\ref{lem:points-v-odd}
and \ref{lem:points-v-3}, except for the case of $b \equiv 16 \bmod 64$ when
$\ord_{2}(x(P))>0$ (when $P$ has singular reduction).

When $b \equiv 16 \bmod 64$ and $\ord_{2}(x(P))>0$, we will use
a minimal model. From Step~11 in the proof of Lemma~\ref{lem:curves-v-2}, if
$P=(x(P),y(P)) \in E_{b} \left( \bbQ_{2} \right)$, then $Q=(x(P)/4,y(P)/8-1/2)
\in E_{b,min} \left( \bbQ_{2} \right)$ defined by $y^{2}+y=x^{3}+(b-16)/64$.
$Q$ has nonsingular reduction and so
\begin{align*}
\widehat{\lambda}_{2}(P)=\widehat{\lambda}_{2}(Q)
&= \frac{1}{2} \log \max \left\{ 1, |x(P)/4|_{2} \right\}
   - \frac{1}{12}\log \left| \Delta \left( E_{b,min} \right) \right|_{2} \\
&= \frac{1}{2} \log \max \left\{ 1, |x(P)|_{2} \right\}
   - \frac{1}{12}\log \left| \Delta \left( E_{b} \right) \right|_{2}-\log(2),
\end{align*}
since $\ord_{2}(x(P))>0$ implies $\ord_{2}(x(P)) \geq 2$ (since $\ord_{2}(b)=4$)
and $\Delta \left( E_{b,min} \right)=\Delta \left( E_{b} \right)/2^{12}$.
\end{proof}

\subsection{Contributions from $p=2$ and $3$}
\label{ssec:2-and-3}

It will be useful for the proofs of our theorems to collect the information
required for the contributions from $p=2$ and $p=3$. We do so here in
Tables~\ref{table:quant-2} and \ref{table:quant-3}.

Below, $C_{3,2}$ and $C_{3,3}$ are the reciprocals of the exponentials of the
quantities from Lemmas~\ref{lem:points-v-3}(b) and \ref{lem:points-v-2}(b) above.

For a point $P \in E_{b}(\bbQ)$, $C_{4,p}=p^{\max \left( 0, \ord_{p}(x(P)) \right)}$.
Similarly, $C_{5,p}=p^{\max \left( 0, 2\ord_{p}(y(P)) \right)}$ and
$C_{6,p}=p^{\max \left( 0, -\ord_{p}(x(2P)) \right)}$. Our computations showed
that $C_{6,3}=1$ is possible for each entry in our table below, so we leave it out
and write $C_{6}$ for $C_{6,2}$.

Letting $\ord_{2}(b)=k$ and $\ord_{3}(b)=\ell$, we put $C_{7,2}=C_{3,2}/2^{k/6}$
and $C_{7,3}=C_{3,3}/3^{\ell/6}$.

We define $C_{3}=C_{3,2}C_{3,3}$, $C_{4}=C_{4,2}C_{4,3}$, $C_{5}=C_{5,2}C_{5,3}$
and $C_{7}=C_{7,2}C_{7,3}$.

In the tables below, the values of $C_{3,p}$, $C_{4,p}$ and $C_{5,p}$ are
expressed as $\left( v_{1}:v_{2} \right)$, where $v_{1}$ is the minimum
possible value when $x(P)$ has good reduction modulo $p$ and under the conditions
on $b$, while $v_{2}$ is the minimum possible value when $x(P)$ has singular
reduction modulo $p$ and under the conditions on $b$. A ``-'' indicates that
a case that is not possible.

\begin{table}[h]
\centering
\scalebox{0.75}{%
\hspace{-22.0mm}\begin{tabular}{|c|c|c|c|c|c||c|c|c|c|c|c|}
\hline
$b$           &   $C_{3,2}$                & $C_{4,2}$ & $C_{5,2}$ & $C_{6,2}$ &          $C_{7,2}$                &         $b$      & $C_{3,2}$ & $C_{4,2}$ & $C_{5,2}$ & $C_{6,2}$ &          $C_{7,2}$          \\ \hline
$ 1 \bmod  8$ & $\left( 1:2^{1/3} \right)$ & $(1: 2)$  & $( 1: 1)$ & $(16: 1)$ & $\left(        1:2^{1/3} \right)$ &    $32 \bmod 64$ &  $(1:-)$  &  $(1:-)$  & $(1:-)$   & $( 4:-)$  & $\left( 2^{-5/6}:- \right)$ \\ \hline
$ 8 \bmod 16$ & $\left( 1:2^{1/2} \right)$ & $(1: 2)$  & $( 1:16)$ & $( 4: 4)$ & $\left( 2^{-1/2}:      1 \right)$ &    $48 \bmod 64$ &  $(1:-)$  &  $(1:-)$  & $(1:-)$   & $( 4:-)$  & $\left( 2^{-2/3}:- \right)$ \\ \hline
$12 \bmod 16$ & $\left( 1:2^{1/2} \right)$ & $(1: 2)$  & $( 1: 4)$ & $( 4: 1)$ & $\left( 2^{-1/3}:2^{1/6} \right)$ &    $20 \bmod 32$ &  $(1:-)$  &  $(1:-)$  & $(1:-)$   & $( 4:-)$  & $\left( 2^{-1/3}:- \right)$ \\ \hline
$ 4 \bmod 32$ & $\left( 1:2^{2/3} \right)$ & $(1: 4)$  & $( 1: 4)$ & $( 4: 1)$ & $\left( 2^{-1/3}:2^{1/3} \right)$ &     $2 \bmod  4$ &  $(1:-)$  &  $(1:-)$  & $(1:-)$   & $( 4:-)$  & $\left( 2^{-1/6}:- \right)$ \\ \hline
$16 \bmod 64$ & $\left( 1:      2 \right)$ & $(1: 4)$  & $( 1:16)$ & $( 4: 1)$ & $\left( 2^{-2/3}:2^{1/3} \right)$ & $3,5,7 \bmod  8$ &  $(1:-)$  &  $(1:-)$  & $(4:-)$   & $(16:-)$  & $\left(        1:- \right)$ \\ \hline
\end{tabular}%
}
\caption{Quantities for $p=2$}
\label{table:quant-2}
\end{table}

%
\begin{table}[h]
\centering
\scalebox{0.75}{%
\hspace{-16.0mm}\begin{tabular}{|c|c|c|c|c||c|c|c|c|c|c|c|}
\hline
$b$                &         $C_{3,3}$          & $C_{4,3}$ &  $C_{5,3}$ &          $C_{7,3}$                &           $b$       & $C_{3,3}$ & $C_{4,3}$ & $C_{5,3}$ &         $C_{7,3}$           \\ \hline
$ 1,8   \bmod   9$ & $\left( 1:3^{1/4} \right)$ & $(1: 1)$  & $(1,   9)$ & $\left(        1:3^{1/4} \right)$ & $243,486 \bmod 729$ &  $(1:-)$  &  $(1:-)$  &  $(1:-)$  & $\left( 3^{-5/6}:- \right)$ \\ \hline
$ 9     \bmod  27$ & $\left( 1:3^{1/3} \right)$ & $(1: 3)$  & $(1,   9)$ & $\left( 3^{-1/3}:      1 \right)$ &     $162 \bmod 243$ &  $(1,-)$  &  $(1:-)$  &  $(1:-)$  & $\left( 3^{-2/3}:- \right)$ \\ \hline
$54,108 \bmod 243$ & $\left( 1:3^{2/3} \right)$ & $(1: 3)$  & $(1,  81)$ & $\left( 3^{-1/2}:3^{1/6} \right)$ & $135,189 \bmod 243$ &  $(1,-)$  &  $(1:-)$  &  $(1:-)$  & $\left( 3^{-1/2}:- \right)$ \\ \hline
$81     \bmod 243$ & $\left( 1:3^{2/3} \right)$ & $(1: 9)$  & $(1,  81)$ & $\left( 3^{-2/3}:      1 \right)$ & $     18 \bmod  27$ &  $(1,-)$  &  $(1:-)$  &  $(1:-)$  & $\left( 3^{-1/3}:- \right)$ \\ \hline
$27,216 \bmod 243$ & $\left( 1:3^{3/4} \right)$ & $(1: 3)$  & $(1, 729)$ & $\left( 3^{-1/2}:3^{1/4} \right)$ & $    3,6 \bmod   9$ &  $(1,-)$  &  $(1:-)$  &  $(1:-)$  & $\left( 3^{-1/6}:- \right)$ \\ \hline
                   &                            &           &            &                                   & $2,4,5,7 \bmod   9$ &  $(1,-)$  &  $(1:-)$  &  $(1:-)$  & $\left(        1:- \right)$ \\ \hline
\end{tabular}%
}
\caption{Quantities for $p=3$}
\label{table:quant-3}
\end{table}

The values of $C_{4,p}$, $C_{5,p}$ and $C_{6,2}$ are obtained by computation.
Using PARI, we calculate these values for each possibility modulo $p^{6}$ of
$x(P)=\alpha/\delta^{2}$, where $\alpha$ and $\delta$ are relatively prime
integers, and $b$ modulo $p^{6}$.

\subsection{Global minimal Weierstrass equation for $E_{b}/\bbQ$}

Putting together the information we obtained from Tate's algorithm in the above
three subsections we obtain the following result.

\begin{lemma}
\label{lem:minimal}
Let $b_{1}$ be the sixth-power-free part of $b$. If $b_{1} \equiv 16 \bmod 64$,
then a global minimal Weierstrass equation for $E_{b}/\bbQ$ is
\[
y^{2}+y = x^{3}+\left(b_{1}-16 \right)/64.
\]
Otherwise, a global minimal Weierstrass equation for $E_{b}/\bbQ$ is
\[
y^{2}= x^{3}+b_{1}.
\]
\end{lemma}

\section{Proof of Theorem~\ref{thm:lang}}

\subsection{Proof of part~(a) ($c_{p}=1$)}

We compute the canonical height by summing local heights.

From Lemma~\ref{lem:curves-v-odd} and our hypotheses, $P$ has nonsingular
reduction for all $P \in E_{b} \left( \bbQ_{v} \right)$ and all primes $q_{v}>3$.
Hence we can apply Lemma~\ref{lem:points-v-odd}(b) for these primes. Combining
this with Lemmas~\ref{lem:points-v-3}(b) and \ref{lem:points-v-2}(b) gives the
inequality
\begin{equation}
\label{eq:nonarch-sum1}
\sum_{v \neq \infty} \widehat{\lambda}_{v}(P)
\geq - \log \left( C_{3} \right) + \frac{1}{12} \log \left| \Delta \left( E_{b} \right) \right|.
\end{equation}

\vspace*{1.0mm}

{\sc Case} $b<0$. Adding \eqref{eq:nonarch-sum1} to the lower bound obtained from \eqref{eq:archim-hgt-b-neg1}
for $\widehat{\lambda}_{\infty}(P)$, we have
\begin{equation}
\label{eq:lang1-neg}
\widehat{h}(P) > \frac{1}{6} \log |b| - \log \left( C_{3} \right) + \frac{1}{4}\log(3).
\end{equation}

From Tables~\ref{table:quant-2} and \ref{table:quant-3}, we see that the
minimum value of $C_{3}$ is $2 \cdot 3^{3/4}$, which occurs when $b \equiv 16 \bmod 64$
and $b \equiv 27, 216 \bmod 243$.

\vspace*{1.0mm}

{\sc Case} $b>0$. Adding \eqref{eq:nonarch-sum1} to the lower bound from \eqref{eq:archim-hgt-b-pos1}
for $\widehat{\lambda}_{\infty}(P)$, we obtain
\begin{equation}
\label{eq:lang1-pos}
\widehat{h}(P) > \frac{1}{6} \log |b| - \log \left( C_{3} \right) + \frac{1}{3}\log(2)
-0.006.
\end{equation}

\subsection{Proof of part~(b) ($c_{p}|4$)}

Again, we compute the canonical height by summing local heights.

From Lemma~\ref{lem:curves-v-odd} and our hypotheses, $[2]P$ has nonsingular
reduction for all $P \in E_{b} \left( \bbQ_{v} \right)$ and all primes $q_{v}>3$.
Hence we can apply Lemma~\ref{lem:points-v-odd}(b) for these primes. Combining
this with Lemmas~\ref{lem:points-v-3}(b) and \ref{lem:points-v-2}(b), and writing
$x([2]P)=\alpha/\delta^{2}$ as a fraction in lowest terms with $\delta>0$, gives
\begin{eqnarray}
\label{eq:nonarch-sum2}
\sum_{v \neq \infty} \widehat{\lambda}_{v}\left( [2](P) \right)
& \geq & \log (\delta) -\log \left( C_{3}' \right) + \frac{1}{12} \log \left| \Delta \left( E_{b} \right) \right| \nonumber \\
& \geq & \frac{1}{2}\log \left( C_{6} \right) -\log \left( C_{3}' \right) + \frac{1}{12} \log \left| \Delta \left( E_{b} \right) \right|,
\end{eqnarray}
where $C_{3}'$ is the value of $C_{3}$ for $[2]P$ (not $P$). These values can
be different since $c_{3}=2$ for $b \equiv 1, 8 \bmod 9$ and $b \equiv 27,216 \bmod 243$,
so all points have nonsingular reduction, and $c_{2}=2$ for $b \equiv 8,12 \bmod 16$,
so again all points have nonsingular reduction.

Note the worst cases occur for $b \equiv 54,81,108 \bmod 243$ and $b \equiv 16 \bmod 64$,
when $C_{3}'/C_{6}^{1/2}=2 \cdot 3^{2/3}$.

\vspace*{1.0mm}

{\sc Case} $b<0$. Adding \eqref{eq:nonarch-sum2} to the lower bound obtained from \eqref{eq:archim-hgt-b-neg1}
for $\widehat{\lambda}_{\infty}\left( [2](P) \right)$ and using $\widehat{h}\left( [2](P) \right)=4\widehat{h}(P)$,
we get
\begin{equation}
\label{eq:lang2-neg}
\widehat{h}(P)
> \frac{1}{24} \log |b| + \frac{1}{8}\log \left( C_{6} \right) - \frac{1}{4}\log \left( C_{3}' \right) + \frac{1}{16}\log(3).
\end{equation}

In the worst cases, $C_{3}'^{1/4}/ \left( C_{6}^{1/8} 3^{1/16} \right)=2^{1/4} \cdot 3^{5/48}$,
Theorem~\ref{thm:lang}(b) immediately follows in this case.

\vspace*{1.0mm}

{\sc Case} $b>0$. Adding \eqref{eq:nonarch-sum2} to the lower bound from \eqref{eq:archim-hgt-b-pos1}
for $\widehat{\lambda}_{\infty}\left( [2](P) \right)$ and using $\widehat{h}\left( [2](P) \right)=4\widehat{h}(P)$,
we have
\begin{equation}
\label{eq:lang2-pos}
\widehat{h}(P)
> \frac{1}{24} \log |b| + \frac{1}{8}\log \left( C_{6} \right) - \frac{1}{4}\log \left( C_{3}' \right) + \frac{1}{12}\log(2)-0.002.
\end{equation}

Here in the worst cases, $C_{3}'^{1/4} / \left( C_{6}^{1/8} 3^{1/16} \right) =2^{1/6} \cdot 3^{1/6}$,
completing the proof of Theorem~\ref{thm:lang}(b).

\subsection{Proof of part~(d) ($c_{p}|12$)}

We prove part~(d) first as we can then use simplified versions of some of the
statements here in the proof of part~(c).

Write $b=2^{k}3^{\ell}q_{2}^{2}q_{3}^{3}q_{4}^{4}q$ where $q_{2}$ is the
product of all distinct primes, $p \geq 5$, with $\ord_{p}(b)=2$, $\ord_{p}(x(P))>0$
and $b/p^{2}$ a quadratic residue modulo $p$; $q_{3}$ is the product of all
distinct primes, $p \geq 5$, with $\ord_{p}(b)=3$, $\ord_{p}(x(P))>0$ and
$b/p^{3}$ a cubic residue modulo $p$; $q_{4}$ is the product of all distinct
primes, $p \geq 5$, with $\ord_{p}(b)=4$, $\ord_{p}(x(P))>0$ and $b/p^{4}$ a
quadratic residue modulo $p$; and $q$ the remaining divisors of $b$ with
$\gcd \left( 6q_{2}q_{3}q_{4}, q \right)=1$. We put $Q_{2}=q_{2}q_{4}^{2}$.

Notice that if a prime $p$ is a divisor of $q_{4}$ and $\ord_{p}(x(P))>0$,
then, in fact, $\ord_{p}(x(P)) \geq 2$. Otherwise if $p$ is a prime dividing
$q_{4}$ with $\ord_{p}(x(P))=1$, then $\ord_{p} \left( x(P)^{3}+b \right)=3$,
but it must be even (since it equals $\ord_{p} \left( y(P)^{2} \right)$).
Similarly, if $k \geq 4$ or $\ell \geq 4$,
then $\ord_{2}(x(P)) \geq 2$ or $\ord_{3}(x(P)) \geq 2$, respectively.

Writing $x(P)=\alpha/\delta^{2}$ with $\alpha$ and $\delta>0$ relatively prime
integers (see, for example, \cite[\S~III.2]{Silv-Tate}), we have
\begin{equation}
\label{eq:alpha-factor-d}
\left( C_{4}q_{2}'q_{3}'q_{4}'^{2} \right) | \alpha,
\end{equation}
where $C_{4}$ is as above, $q_{2}'|q_{2}$, $q_{3}'|q_{3}$ and $q_{4}'|q_{4}$.
We put $Q_{2}'=q_{2}'q_{4}'^{2}$, noting that $Q_{2}'|Q_{2}$. So we can write
$\alpha=C_{4}Q_{2}'q_{3}'q'$ for an integer $q'$.

We write $x(P)=c|b|^{1/3}$ for $c \geq -1$ and combining this with
\eqref{eq:alpha-factor-d}, we find that
$C_{4}Q_{2}'q_{3}' \leq c|b|^{1/3} \delta^{2}$. That is,
\begin{equation}
\label{eq:q2q3q4-ub2}
Q_{2}'^{2}q_{3}'^{2} \leq c^{2} |b|^{2/3}\delta^{4}/C_{4}^{2}.
\end{equation}

Since
$x(P)^{3}+b= \left( C_{4}q'Q_{2}'q_{3}'/\delta^{2} \right)^{3} + 2^{k}3^{\ell}qQ_{2}^{2}q_{3}^{3}$
is a perfect square, so is
$$
q_{3}' \left( C_{4}^{3}q'^{3}Q_{2}'+2^{k}3^{\ell}q\left( \frac{Q_{2}q_{3}}{Q_{2}'q_{3}'} \right)^{2}\left( q_{3}/q_{3}' \right) \delta^{6} \right)
=q_{3}' Q'.
$$

Since $\gcd \left( 6, q_{3}' \right)=1$, it must be the case that $q_{3}'$
divides $Q'/C_{5}$. Thus
$$
q_{3}' \leq
\left( C_{4}^{3}q'^{3}Q_{2}'+2^{k}3^{\ell}q\left( \frac{Q_{2}q_{3}}{Q_{2}'q_{3}'} \right)^{2}\left( q_{3}/q_{3}' \right) \delta^{6} \right)/C_{5}.
$$

Substituting $q'=c|b|^{1/3}\delta^{2}/\left( C_{4}Q_{2}'q_{3}' \right)$ and our
expression for $b$ into this upper bound for $q_{3}'$, we have
\begin{equation}
\label{eq:q3-ub}
q_{3}' \leq 2^{k}3^{\ell}|q|\delta^{6} \left( c^{3}+\sgn(b) \right)/C_{5}.
\end{equation}

Combining \eqref{eq:q3-ub} with our expression for $b$ to eliminate $q$, we
obtain
\begin{equation}
\label{eq:q2q3q4-ub3}
Q_{2}'^{2}q_{3}'^{4} \leq \frac{\delta^{6} \left( c^{3}+\sgn(b) \right)}{C_{5}}|b|.
\end{equation}

So, from \eqref{eq:q2q3q4-ub2} and \eqref{eq:q2q3q4-ub3}, we have
\begin{equation}
\label{eq:q2q3q4-ub4}
\left( Q_{2}'^{2}q_{3}'^{3} \right)^{2}
\leq \frac{c^{2}\delta^{10} \left( c^{3}+\sgn(b) \right)}{C_{4}^{2}C_{5}}|b|^{5/3}.
\end{equation}

From Lemmas~\ref{lem:points-v-odd}(b), \ref{lem:points-v-3}(b) and
\ref{lem:points-v-2}(b), along with \eqref{eq:alpha-factor-d} and
\eqref{eq:q2q3q4-ub4}, we obtain
\begin{eqnarray}
\label{eq:nonarch-sum4}
\sum_{v \neq \infty} \widehat{\lambda}_{v}(P)
& \geq & \log (\delta) - \frac{\log \left( Q_{2}'^{2}q_{3}'^{3} \right)}{6}
         - \log \left( C_{4} \right)
         + \frac{\log \left| \Delta \left( E_{b} \right) \right|}{12} \nonumber \\
& \geq & \log (\delta) - \frac{\log \left( c^{2}\delta^{10} \left( c^{3}+\sgn(b) \right)|b|^{5/3} / \left( C_{4}^{2}C_{5} \right) \right)}{12} \nonumber \\
&      & - \log \left( C_{3} \right)
         + \frac{\log \left| \Delta \left( E_{b} \right) \right|}{12}.
\end{eqnarray}

\vspace*{1.0mm}

{\sc Case} $b<0$. Here we have $x(P)=c|b|^{1/3}$ for $c \geq 1$.

For $c>1$, we combine \eqref{eq:nonarch-sum4} with \eqref{eq:archim-hgt-b-neg3b}
in Lemma~\ref{lem:arch-b-neg} obtaining
\begin{eqnarray}
\label{eq:lang4-neg}
\widehat{h}(P)
&   >  & \frac{\log |b|}{6} + \frac{1}{12} \log \left( c^{5}-c^{2} \right) + 0.1895
         + \log (\delta) \nonumber \\
&      & - \frac{\log \left( c^{2}\delta^{10} \left( c^{3}-1 \right)|b|^{5/3} / \left( C_{4}^{2}C_{5} \right) \right)}{12}
         - \log \left( C_{3} \right)
         + \frac{\log \left| \Delta \left( E_{b} \right) \right|}{12} \nonumber \\
& \geq & \frac{\log |b|}{36} + \frac{\log \left( C_{4}^{2}C_{5}/C_{3}^{12} \right)}{12} + 0.1895,
\end{eqnarray}
since $\delta \geq 1$.

From Tables~\ref{table:quant-2} and \ref{table:quant-3}, the minimum value of
$C_{4}^{2}C_{5}/C_{3}^{12}$ is $2^{-4} \cdot 3^{-2}$, which can occur for
$b \equiv 54,108 \bmod 243$ and $b \equiv 16 \bmod 64$.

Note that in Lemma~\ref{lem:arch-b-neg}(b), we exclude $c=1$. However, this
is a torsion point, which is excluded from our results (see the argument in
the next section using \cite{F}).

\vspace*{1.0mm}

{\sc Case} $b>0$. Here we have $x(P)=cb^{1/3}$ for $c \geq -1$.

The argument is identical to that for $b<0$, except that we use
\eqref{eq:archim-hgt-b-pos3b} in Lemma~\ref{lem:arch-b-pos}, rather than
\eqref{eq:archim-hgt-b-neg3b} in Lemma~\ref{lem:arch-b-neg}. Thus
\begin{equation}
\label{eq:lang4-pos}
\widehat{h}(P) >
\frac{\log |b|}{36} + \frac{\log \left( C_{4}^{2}C_{5}/C_{3}^{12} \right)}{12} + 0.188.
\end{equation}

In Lemma~\ref{lem:arch-b-pos}(b), we exclude $c=-1$ and $c=0$. But, as above,
these are torsion points and are not under consideration here.

Hence the theorem holds for $b>0$ too.

\subsection{Proof of part~(c) ($c_{p}|3$)}

We proceed as in the proof of part~(d), using the notation there too, except
here we have $q_{3}=q_{3}'=1$. Thus
\begin{equation}
\label{eq:nonarch-sum3}
\sum_{v \neq \infty} \widehat{\lambda}_{v}(P)
\geq \log (\delta) - \frac{\log \left( Q_{2}' \right)}{3}
- \log \left( C_{3} \right)
+ \frac{\log \left| \Delta \left( E_{b} \right) \right|}{12}.
\end{equation}
and
\begin{equation}
\label{eq:q2q4-ub2}
Q_{2}'^{2} \leq c^{2} |b|^{2/3}\delta^{4}/C_{4}^{2}.
\end{equation}

\vspace*{1.0mm}

{\sc Case} $b<0$. Note that here $c \geq 1$.

We combine \eqref{eq:nonarch-sum3} with \eqref{eq:archim-hgt-b-neg3a} in
Lemma~\ref{lem:arch-b-neg}(b) to obtain
\begin{eqnarray}
\label{eq:lang3-bNeg-a}
\widehat{h}(P)
& > & \frac{\log |b|}{6} + \frac{\log(c)}{3} + \frac{\log(108)}{18} \\
&   & - 0.004 + \log (\delta) - \frac{\log \left( Q_{2}' \right)}{3}
  - \log \left( C_{3} \right). \nonumber
\end{eqnarray}

Now we apply $\delta \geq 1$ and \eqref{eq:q2q4-ub2} to \eqref{eq:lang3-bNeg-a},
obtaining
\begin{eqnarray}
\label{eq:lang3-neg}
\widehat{h}(P)
&   >  & \frac{\log |b|}{6} + \frac{\log(c)}{3} + \frac{\log(108)}{18} - 0.004
         + \log (\delta) - \log \left( C_{3} \right) \nonumber \\
&      & - \frac{\log(c)}{3} - \frac{\log|b|}{9} - \frac{2\log (\delta)}{3} + \frac{\log \left( C_{4} \right)}{3} \nonumber \\
& \geq & \frac{1}{18}\log |b| + \frac{\log (108)}{18}
+ \frac{\log\left( C_{4}/C_{3}^{3} \right)}{3} - 0.004.
\end{eqnarray}

Again, from Tables~\ref{table:quant-2} and \ref{table:quant-3}, the minimum
value of $C_{4}/C_{3}^{3}$ is $2^{-4} \cdot 3^{-9/2}$,
which can occur for $b \equiv 27, 216 \bmod 243$ and $b \equiv 16 \bmod 64$.

\vspace*{1.0mm}

{\sc Case} $b>0$. We proceed in the same way as for $b<0$, except using
\eqref{eq:archim-hgt-b-pos3a} in Lemma~\ref{lem:arch-b-pos}(b), to obtain
\begin{equation}
\label{eq:lang3-pos}
\widehat{h}(P)
> \frac{1}{18}\log |b| + \frac{\log (27)}{12}
  + \frac{\log \left( C_{4}/C_{3}^{3} \right)}{3} - 0.004.
\end{equation}

The minimum value of $27C_{4}^{4}/C_{3}^{12}$ is
$2^{-4} \cdot 3^{-2}$, which can occur for $b \equiv 27, 216 \bmod 243$ and
$b \equiv 16 \bmod 64$.

We also record here the analogue of \eqref{eq:lang3-bNeg-a} which can be useful
in many specific cases.
\begin{equation}
\label{eq:lang3-bNeg-b}
\widehat{h}(P)
> \frac{\log |b|}{6} + \frac{\log(c)}{3}
  + \frac{\log(27)}{27} - 0.004
  + \log (\delta) - \frac{\log \left( Q_{2}' \right)}{3}
  - \log \left( C_{3} \right).
\end{equation}

As in the proof of part~(d), where appropriate, the points with $c=-1$, $0$
or $1$, correspond to torsion points and are not considered here.
Hence part~(c) of Theorem~\ref{thm:lang} holds too.

\section{Proof of Theorem~\ref{thm:hgt-diff}}

As in the previous section, write $b=2^{k}3^{\ell}q_{2}^{2}q_{3}^{3}q_{4}^{4}q$
and $x(P)=C_{4}q_{2}'q_{3}'q_{4}^{2}q'$.
From Lemmas~\ref{lem:points-v-odd}(b), \ref{lem:points-v-3}(b) and
\ref{lem:points-v-2}(b), along with the definitions of $C_{3}$ and $C_{7}$ in
Subsection~\ref{ssec:2-and-3}, we get
\begin{align}
\label{eq:non-arch}
0 &\leq \sum_{v \neq \infty} \left(
         \frac{1}{2} \log \max \left\{ 1, |x(P)|_{v} \right\}
         - \frac{1}{12} \log \left| \Delta \left( E_{b} \right) \right|_{v}
         - \widehat{\lambda}_{v}(P) \right) , \nonumber \\
  &=    \frac{1}{3}\log \left| q_{2}' \right| + \frac{1}{2}\log \left| q_{3}' \right|
        + \frac{2}{3}\log \left| q_{4}' \right| - \frac{k}{6}\log(2) - \frac{\ell}{6}\log(3)+\log \left( C_{3} \right) \\
  &\leq \frac{1}{6}\log |b| + \log \left( C_{7} \right), \nonumber
\end{align}
with the upper bound achieved when for every prime $p>3$ that divides $b$,
we are in one of the first three cases of \eqref{eq:vadic-not-3-hgt}. That is
$|q'|=|q|=1$, $q_{2}'=q_{2}$, $q_{3}'=q_{3}$ and $q_{4}'=q_{4}$.

If $b<0$, then from Lemma~\ref{lem:arch-b-neg}(c) and \eqref{eq:non-arch}
\begin{equation}
\label{eq:diff-neg}
-\frac{\log(3)}{4}-0.005
< \frac{1}{2}h(P) - \widehat{h}(P)
< \frac{1}{6}\log|b| + \log \left( C_{7} \right).
\end{equation}

Note from Tables~\ref{table:quant-2} and \ref{table:quant-3} that the maximum
value of $C_{7}=2^{1/3} \cdot 3^{1/4}$, which can occur for $b \equiv 1 \bmod 8$,
$4 \bmod 32$ or $16 \bmod 64$ and $b \equiv 1,8 \bmod 9$, or $27,216 \bmod 243$.

Now suppose that $b \geq 2$, then from Lemma~\ref{lem:arch-b-pos}(c) and
\eqref{eq:non-arch}
\begin{eqnarray}
\label{eq:diff-pos}
&   & -\frac{1}{6}\log |b| - \frac{\log(2)}{3} - 0.007-0.076b^{-1/3} \nonumber \\
& < & \frac{1}{2}h(P) - \widehat{h}(P)
      < \frac{1}{6}\log|b| + \log \left( C_{7} \right) + 0.004.
\end{eqnarray}

For $b=1$, $E_{b}(\bbQ)$ consists only of torsion points, which we consider
next for all $b$.

From \cite{F} (see also \cite[Proposition~6.31]{S-Z}), the torsion group
of $E_{b}(\bbQ)$ is isomorphic to:\\
$\bullet$ $\bbZ/6\bbZ$, if $b=1$
(the torsion points are $(2,\pm 3)$, $(0, \pm 1)$, $(-1,0)$, $O$) \\
$\bullet$ $\bbZ/3\bbZ$, if $b=b_{1}^{2} \neq 1$ or if $b=-432$
(the torsion points are $\left(0, \pm b_{1} \right)$, $O$ in the former case
and $(12,\pm 36)$, $O$ in the latter) \\
$\bullet$ $\bbZ/2\bbZ$, if $b=b_{1}^{3} \neq 1$
(the torsion points are $\left( -b_{1},0 \right)$ and $O$) \\
$\bullet$ $\{ O \}$, otherwise.

In the first case, $0 \leq (1/2)h(P)-\widehat{h}(P) \leq \log(2)/2$.

In the second case when $b=b_{1}^{2}$, $h(P)=\widehat{h}(P)=0$.

In the second case when $b=-432$, $0 \leq (1/2)h(P)-\widehat{h}(P) \leq \log(12)/2$.

Lastly, in the third case, $0 \leq (1/2)h(P)-\widehat{h}(P) \leq \log \left( b_{1} \right)/2
=\log(b)/6$.

So in all these cases, our Theorem holds as well.

For the last inequality in Theorem~\ref{thm:hgt-diff}, we observe that for $b=\pm 1$,
$E_{b}(\bbQ)$ contains only the torsion points, so we may assume that $|b| \geq 2$.

For $b \leq -2$, $-(1/6) \log |b|-0.299 < -0.41 < -\log(3)/4-0.005$.

For $b \geq 2$, $-\log(2)/3-0.007-0.076/b^{1/3}=-0.298\ldots$, so this
inequality holds.

\section{Sharpness of Results}
\label{sect:sharp}

For each part of our Theorems, we produce infinite families of pairs of curves
and points on those points demonstrating that the results, without the small
constant ``error terms'', are best possible (excluding Theorem~\ref{thm:lang}(d)
where our examples are within a very small constant of what we believe are the
best possible results).

\subsection{Theorem~\ref{thm:lang}(a)}

\quad \\
{\sc Case} $b<0$. Set
$$
b=-46656b_{1}^{3} - 93312b_{1}^{2} - 62208b_{1} - 2160
\mbox{ and } P= \left( 36b_{1}+24, 108 \right),
$$
where $b_{1}$ is a positive integer and we let it approach $+\infty$. We find
that $x(P) \rightarrow |b|^{1/3}$ and hence the archimedean height approaches
the lower bound in Lemma~\ref{lem:arch-b-neg}(a). Since $b \equiv 16 \bmod 64$
and $b \equiv 27 \bmod 243$, such values of $b$ have the smallest nonarchimedean
height functions at both $2$ and $3$. Furthermore, by our conditions on $b$ in
Theorem~\ref{thm:lang}(a), our points $P$ have nonsingular reduction for the
other primes.

\vspace*{1.0mm}

{\sc Case} $b>0$. Take
$$
b=46656b_{1}^{2}+46656b_{1}+13392 \hspace*{3.0mm} \text{ and } \hspace*{3.0mm}
P=\left( -12, 54 \left( 4b_{1}+2 \right) \right),
$$
where $b_{1}$ is a positive integer and we let it approach $+\infty$.

For such pairs of curves and points, we find that $x(P)/|b|^{1/3} \rightarrow 0$
as $b_{1} \rightarrow +\infty$ and hence the archimedean height approaches the
lower bound in Lemma~\ref{lem:arch-b-pos}(a). As in the case of $b<0$, the
required conditions at each of the primes are satisfied too.

\subsection{Theorem~\ref{thm:lang}(b)}

\quad \\
{\sc Case} $b<0$. Let $b_{1}$ be an odd positive integer and put
$b_{2}=\left[ 12b_{1}^{3}/ \left( 3+2\sqrt{3} \right) \right]$ where $[z]$ is
the nearest integer to $z$. Put
$$
b=-432\left( 12b_{1}^{3}-b_{2} \right)^{3} \left( 3b_{2}-4b_{1}^{3} \right)
\mbox{ and }
P= \left( 24b_{1} \left( 12b_{1}^{3}-b_{2} \right), 36\left( 12b_{1}^{3}-b_{2} \right)^{2} \right).
$$

Suppose that $b_{1}$ and $b_{2}$ are relatively prime and that $2^{6} \nmid b$
and $3^{6} \nmid b$.

For such pairs, $x\left( [2](P) \right)=48b_{1}b_{2}$ and $x\left( [2](P) \right)^{3} \approx 191102976b_{1}^{12}/(3+2\sqrt{3})^{3}$.
We also find that
$b=-191102976b_{1}^{12}/(3+2\sqrt{3})^{3}+O\left(b_{1}^{9} \right)$. Therefore,
as $b_{1} \rightarrow +\infty$, $x\left( [2](P) \right) \rightarrow |b|^{1/3}$ and the lower
bound for the archimedean height is sharp.

\vspace*{1.0mm}

{\sc Case} $b>0$. Let $b_{1}$ be a positive integer, put
$$
b=432\left( 162b_{1}+31 \right) \left( 6b_{1}+1 \right)^{3}
\mbox{ and }
P = \left( -72b_{1}-12, 108 \left( 6b_{1}+1 \right)^{2} \right).
$$

For such pairs, $x\left( [2](P) \right)=144b_{1}+28$, so $x\left( [2](P) \right)/|b|^{1/3} \rightarrow 0$ as
$b \rightarrow +\infty$ and the lower bound for the archimedean height in
Lemma~\ref{lem:arch-b-pos}(a) is sharp as $b_{1} \rightarrow +\infty$. As above,
the desired conditions on all the primes are satisfied too (note
$b \equiv 16 \bmod 64$ and $b \equiv 27 \bmod 243$).

\subsection{Theorem~\ref{thm:lang}(c)}

\quad \\
{\sc Case} $b<0$. Let $b_{1}$ be a positive integer and put $b_{2}=36b_{1}^{2}+36b_{1}+11$,
$$
b=-432 \left( b_{2}+6 \right) b_{2}^{2}
\hspace*{1.0mm} \mbox{ and } \hspace*{1.0mm}
P= \left( 12b_{2}, 324 \left( 2b_{1}+1 \right) b_{2} \right).
$$

Here $x \left( \left[ 2^{n} \right](P) \right) \rightarrow 4^{1/3}|b|^{1/3}$ and
$z \left( \left[ 2^{n} \right](P) \right) \rightarrow 3$ as $b_{1} \rightarrow +\infty$
for $n \geq 0$. Therefore
$\widehat{\lambda}_{\infty}(P) \rightarrow (1/6)\log|b|+(1/6)\log(12)-(1/12) \log \left| \Delta_{b} \right|$.

Note that $b \rightarrow -432b_{2}^{3}$, so the sum of the nonarchimedean
heights is
$-(1/9)\log |b| - (5/9)\log(2) - (5/12)\log(3)+(1/12)\log \left| \Delta_{b} \right|$.
Combining this with the above, we find that $\widehat{h}(P) \rightarrow
(1/18)\log|b|-(2/9)\log(2)-(1/4)\log(3)$ from above as $b_{1} \rightarrow +\infty$.

\vspace*{1.0mm}

{\sc Case} $b>0$. Let $b_{1}$ be a positive integer, put $b_{2}=24b_{1}^{2}+24b_{1}+5$,
$$
b=216 \left( b_{2}+9 \right) b_{2}^{2}
\hspace*{1.0mm} \mbox{ and } \hspace*{1.0mm}
P=\left( 12b_{2}, 324 \left( 2b_{1}+1 \right) b_{2} \right).
$$

Here $x(P) \rightarrow 2b^{1/3}$ as $b_{1} \rightarrow +\infty$. We translate
the point and have $x\left( P' \right) \rightarrow 4b^{1/3}$. So $x\left( P' \right)^{4}z\left( P' \right) \rightarrow 72b^{4/3}$.
Furthermore $x \left( \left[ 2^{n} \right]\left( P' \right) \right) \rightarrow 2b^{1/3}$ and
$z \left( \left[ 2^{n} \right]\left( P' \right) \right) \rightarrow 1/2$ for $n \geq 1$. Therefore
$\widehat{\lambda}_{\infty}\left( P' \right) \rightarrow (1/6)\log(b)+(1/8)\log(72)-(1/24)\log(2)
-(1/12) \log \left| \Delta_{b} \right|$.

Note that $b \rightarrow 216b_{2}^{3}$, so the sum of the nonarchimedean
heights is
$-(1/9)\log |b| - (2/3)\log(2) - (5/12)\log(3)+(1/12)\log \left| \Delta_{b} \right|$.
So in this case, $\widehat{h}(P) \rightarrow (1/18)\log|b|-(1/3)\log(2)-(1/6)\log(3)$
as $b_{1} \rightarrow +\infty$.

\subsection{Theorem~\ref{thm:lang}(d)}

Here we produce families where the constants are slightly larger than in the theorem.

\vspace*{1.0mm}

{\sc Case} $b<0$. Let $k$ be a positive integer and put $b_{1}=54k-1$, $b_{2}=720k-1$
and $b_{3}=942k-1$. Note that they are relatively prime and none of them are divisible
by $2$ or $3$. Further, assume that $b_{1}$ and $b_{3}$ are square-free and that $b_{2}$
is cube-free. Let $b=-432b_{1}b_{2}^{2}b_{3}^{3}$ and $P= \left( 12b_{2}b_{3}, 36b_{2}b_{3}^{2} \right)$.

As $k$ increases, $x(P)/|b|^{1/3}$ approaches $(160/3)^{1/3}=3.764\ldots$.
Hence $\widehat{\lambda}_{\infty}(P) \rightarrow (1/6)\log|b|+0.74341680776086\ldots-(1/12)\log \left| \Delta_{b} \right|$.

The sum of the nonarchimedean heights is
$-(1/3)\log \left( b_{2} \right) - (1/2)\log \left( b_{3} \right) - \log(2) - (2/3)\log(3)+(1/12)\log \left| \Delta_{b} \right|$.
Now
\[
b_{2}^{6} \rightarrow \frac{2^{8} \cdot 5^{4}}{3 \cdot 157^{3}}|b|
\hspace*{3.0mm} \text{ and } \hspace*{3.0mm}
b_{3}^{6} \rightarrow \frac{157^{3}}{2^{10} \cdot 3^{7} \cdot 5^{2}}|b|,
\]
so we find that
\[
\widehat{h}(P) \rightarrow \frac{1}{36}\log|b|-0.221457178\ldots
\hspace*{3.0mm} \text{as $k \rightarrow +\infty$}.
\]

Here the constant is approximately $2 \cdot 10^{-7}$ larger than the conjectured
constant. The actual value of $b_{1}$ required to obtain the constant in the
conjecture is smaller than we used here. Here we have $b_{1} \approx (3/40)b_{2}$,
whereas for the conjecture, we require $b_{1} \approx 0.074429578933\ldots b_{2}$.

\vspace*{1.0mm}

{\sc Case} $b>0$. We proceed here just as for $b<0$.

Let $k$ be a positive integer and put $b_{1}=54k-1$, $b_{2}=720k+1$ and $b_{3}=978k+1$.
Note that they are relatively prime and none of them are divisible by $2$ or
$3$. Further, assume that $b_{1}$ and $b_{3}$ are square-free and that $b_{2}$
is cube-free. Let $b=432b_{1}b_{2}^{2}b_{3}^{3}$ and $P= \left( 12b_{2}b_{3}, 36b_{2}b_{3}^{2} \right)$.

With this family of examples, we obtain
\[
\widehat{h}(P) \rightarrow \frac{1}{36}\log|b|-0.22252005826\ldots,
\hspace*{3.0mm} \text{as $k \rightarrow +\infty$}.
\]

As for $b<0$, the constant here is slightly larger than the conjectured constant,
and for the same reason. Here we require $b_{1} \approx 0.085629143\ldots b_{2}$.

\subsection{Theorem~\ref{thm:hgt-diff}}

Silverman (see \cite[Example~2.1]{Silv5}) shows that the coefficients of the
$\log |b|$ terms are best possible.

For the upper bound for $b<0$, we consider
$b=-2^{2} \cdot 3^{3} \cdot 5^{3} b_{1}^{2}$ where
$b_{1}=2160b_{2}^{2}+1350b_{2}+211$ and
$P= \left( 60b_{1}, 1350b_{1} \left( 16b_{2}+5 \right) \right)$,
with the condition that $b_{1}$ be cube-free.

For the upper bound for $b>0$, we consider
$$
b=b_{1}^{2} \hspace*{3.0mm} \text{ and } \hspace*{3.0mm} P=\left( 2b_{1}, 3b_{1} \left( 8b_{2}+15 \right) \right),
$$
where $b_{1}=\left( 6b_{2}+11 \right) \left( 12b_{2}+23 \right)$ and is cube-free.

For the lower bound for $b>0$, we consider
$b= \left( 3b_{1}+1 \right)^{2}+1$ and $P=\left( -1, 3b_{1}+1 \right)$.

For the lower bound for $b<0$, we consider
$b=1-\left( 2b_{1}+1 \right)^{3}$, $P= \left( 2b_{1}+1, 1 \right)$.

In the last inequality of the theorem, $-0.299$ cannot be replaced by anything greater than
$-0.29228\ldots$. Indeed, consider the point $[956](-1,1)$ on $y^{2}=x^{3}+2$
(note that $x\left( [956](-1,1) \right)=0.99818\ldots$). Taking the archimedean
height function evaluated at $x=1$ for $b=2$, we see that $-0.29250\ldots$ is
the smallest possible constant.

\section{Acknowledgements}

The authors would like to express their gratitude to the anonymous referee
for their careful reading of our manuscript, as well as their many useful
suggestions that led to significant improvements in this paper.

\bibliographystyle{amsplain}

\end{document}